\documentclass{amsart}
\usepackage{amsmath}
\usepackage{amsfonts}
\usepackage{graphics}
\usepackage{epsfig}
\usepackage{amssymb}
\usepackage{amscd}

\newtheorem{thm}{Theorem}[section]

\newtheorem{lem}[thm]{Lemma}
\newtheorem{prop}[thm]{Proposition}
\theoremstyle{definition}
\newtheorem{Def}[thm]{Definition}
\newtheorem*{ack}{Acknowledgement}
\theoremstyle{remark}
\newtheorem*{rem}{Remark}

\theoremstyle{definition}
\newtheorem{ex}{Example}[section]

\numberwithin{equation}{section}
\numberwithin{figure}{section}

\def\Hom{{\text{\rm{Hom}}}}
\def\End{{\text{\rm{End}}}}
\def\trace{{\text{\rm{trace}}}}
\def\tr{{\text{\rm{tr}}}}
\def\rchi{{\hbox{\raise1.5pt\hbox{$\chi$}}}}
\def\Aut{{\text{\rm{Aut}}}}

\def\isom{\cong}
\def\tensor{\otimes}
\def\dsum{\oplus}

\def\reg{{\text{\rm{reg}}}}
\def\Ghat{\hat{G}}
\def\lam{\lambda}

\def\fin{{\text{\rm{finite}}}}

\begin{document}
\large
\setcounter{section}{-1}

\title[Non-commutative
matrix integrals and representation varieties]{Non-commutative
matrix integrals and representation varieties of
surface groups in a finite group}
\author[Motohico Mulase]{Motohico Mulase$^1$}  
\address{
Department of Mathematics\\
University of California\\
Davis, CA 95616--8633}
\email{mulase@math.ucdavis.edu}
\author[Josephine Yu]{Josephine T.\ Yu$^2$}  
\address{
Department of Mathematics\\
University of California\\
Berkeley, CA 94720-3840}
\email{jyu@math.berkeley.edu}
\begin{abstract}
A graphical expansion formula for  non-commutative
matrix integrals with values in a finite-dimensional
real or  complex
von Neumann algebra is obtained in terms of  ribbon graphs
and their non-orientable counterpart called
M\"obius graphs. The contribution of 
each graph is an invariant of the topological type of the surface
on which the graph is drawn.
As an example, we calculate the integral on
the  group algebra of a finite group.
We show that
the integral is a generating function of 
the number of homomorphisms from 
the fundamental group of an arbitrary
 closed surface into the finite
group.
The graphical expansion formula yields 
 a new proof of the classical
theorems of Frobenius, Schur and Mednykh on
these numbers. 
\end{abstract}

\thanks{$^1$Research supported   by 
NSF grants DMS-9971371 and DMS-0406077,
and by UC Davis.}
\thanks{$^2$Research supported in part  by NSF grant VIGRE 
DMS-0135345 and UC Davis.}
\maketitle

\allowdisplaybreaks

\tableofcontents

\bigskip

\section{Introduction}
\label{sect:intro}

The purpose of this paper is to establish 
 Feynman diagram expansion formulas
for non-commutative matrix integrals over 
a finite-dimensional real or complex
von Neumann algebra. 
An interesting case is 
 the real or complex group algebra of a
finite group.
 Using the graphical expansion formulas,
we give a new proof of the classical formulas 
for the number of homomorphisms from the fundamental
group of a closed surface
into a finite group,
expressing the number
in terms of  irreducible representations of the finite group.
Indeed, our integrals are generating functions for the
cardinality of the 
\emph{representation variety} of  a surface group
in a finite group. 

The   non-commutative matrix 
integrals of this article
have their origin in \emph{random matrix theory} 
(cf.~\cite{AV, BI, Mehta, TW, VM}), and
  include real symmetric, 
complex hermitian, and quaternionic self-adjoint matrix integrals
as a special case for a \emph{simple} von Neumann algebra. 
Recently a surprising relation between random matrices and random 
permutations was discovered in \cite{BDJ}, and was further studied
from various points of view including representation theory
of symmetric groups (cf.~\cite{BDJ2, BR, BOO, D, J, O, OP1, OP2}).
Our theory exhibits yet
another connection between
matrix-type integrals and representation
theory of finite groups.

Let $A$ be a finite-dimensional complex von Neumann algebra
with the adjoint operation $*:A\rightarrow A$ and a linear
map $\langle \;\rangle:A\rightarrow \mathbb{C}$ called
the \emph{trace}. The algebra $A$ has a positive definite hermitian
inner product defined by 
$$
\langle a,b\rangle = \langle ab^*\rangle
$$
for $a,b\in A$. Let us choose an orthonormal basis 
$\{e_1,\dots,e_N\}$ for $A$ with respect to the hermitian 
form, where $N=\dim A$.  
A \emph{ribbon graph} is a graph with a cyclic order
given at every vertex to incident half-edges.
Recall that every 
ribbon graph
 $\Gamma$ defines a unique closed oriented surface
$S_\Gamma$ on which $\Gamma$ is drawn and gives a 
cell-decomposition. Let $g(\Gamma)$ and $f(\Gamma)$ 
denote the genus of $S_\Gamma$ and the number of $2$-cells,
or faces, of the cell-decomposition, respectively. The Feynman diagram expansion formula
we establish is the following:
\begin{multline}
\label{eq:introexp}
\log \int_{\{x\in A\;|\; x=x^*\}}
\exp\left({-\frac{1}{2}\;\langle x^2\rangle}\right)
\exp\left(\sum_{j=1} ^\infty{\frac{t_j}{j}\;\langle x^j\rangle}\right)
d\mu(x)\\
=
\sum_{\substack{\Gamma  \;\text{ connected}\\
\text{ribbon graph}}}
\frac{1}{|\Aut_R(\Gamma)|}
\; A_{g(\Gamma),f(\Gamma)} ^{or}\;
\prod_j t_j ^{v_j(\Gamma)}\;,
\end{multline}
where
$d\mu(x)$ is a normalized Lebesgue measure on the real vector subspace
of $A$ consisting of self-adjoint elements, $\Aut_R(\Gamma)$ is the
ribbon graph automorphism group,
 and $v_j(\Gamma)$ is the number of 
$j$-valent vertices of the connected ribbon graph $\Gamma$. 
The integral of (\ref{eq:introexp}) is considered as a generating
function of integrals
$$
\int_{\{x\in A\;|\; x=x^*\}}
\exp\left({-\frac{1}{2}\;\langle x^2\rangle}\right)
\prod_j ^\fin \left\langle x^j\right\rangle^{v_j}d\mu(x)
$$
for all finite sequences $(v_1,v_2,v_3,\dots)$ of positive integers.
The contribution of
the graph $\Gamma$ in  (\ref{eq:introexp}) is defined by
$$
A_{g,f} ^{or} = \sum_{\substack{i_1,\dots,i_g;\; j_1,\dots,j_g\\
h_1,\dots,h_{f-1} = 1}} ^N
\left\langle 
e_{i_1}e_{j_1}e_{i_1} ^* e_{j_1} ^*\cdots
e_{i_g}e_{j_g}e_{i_g} ^* e_{j_g} ^*\cdot
e_{h_1}e_{h_1}^* \cdots e_{h_{f-1}} e_{h_{f-1}}^*
\right\rangle\;.
$$
We notice that the graph contribution  
$A_{g(\Gamma),f(\Gamma)} ^{or}$ depends only on the 
\emph{topological type} of the surface $S_\Gamma$, which
is the genus of the surface
and the number of $2$-cells in its cell-decomposition. 
If we apply (\ref{eq:introexp})
 to a simple von Neumann algebra
$A=M(n,\mathbb{C})$, then the formula recovers the well-known
graphical expansion formula for $n\times n$ hermitian matrix
integrals found in many articles, including
\cite{BIZ, Harer-Zagier, Kontsevich, 
Mulase1995, O, OP1, Penner, W1991a}.
The word \emph{non-commutative matrix integral}
in the title is justified because our von Neumann
algebra can take the form $A=B\tensor M(n,\mathbb{C})$ with
another von Neumann algebra $B$. 

For a real von Neumann algebra $A$ with a real valued trace,
our expansion formula is more complicated. Let us 
recall the notion of \emph{M\"obius graph} introduced in 
\cite{MW}. It is the non-orientable counterpart of  ribbon graphs.
A M\"obius graph $\Gamma$ 
defines a unique unoriented surface $S_\Gamma$
and gives a cell-decomposition. Every closed non-orientable 
surface $S$ is obtained by removing $k$ disjoint disks from a sphere
$S^2$ and gluing a \emph{cross-cap} to each hole. The number of 
cross-caps is the \emph{cross-cap genus} of the surface, 
and its Euler characteristic is given by
$\rchi(S) = 2-k$. Every ribbon graph is an orientable 
M\"obius graph, but it has a different automorphism group
reflecting the fact that orientation-reversing map is allowed. 
Now the formula for a real von Neumann algebra is the following:
\begin{equation}
\label{eq:introexpreal}
\begin{split}
\log &\int_{\{x\in A\;|\; x=x^*\}}
\exp\left({-\frac{1}{4}\;\langle x^2\rangle}\right)
\exp\left(\sum_{j=1} ^\infty{\frac{t_j}{2j}\;\langle x^j\rangle}\right)
d\mu(x)\\
&=
\sum_{\substack{\Gamma  \;\text{ connected orientable}\\
\text{M\"obius graph}}}
\frac{1}{|\Aut_M(\Gamma)|}
\; A_{g(\Gamma),f(\Gamma)} ^{or}\;
\prod_j t_j ^{v_j(\Gamma)}\\
&+
\sum_{\substack{\Gamma  \;\text{ connected non-}\\
\text{orientable M\"obius graph}}}
\frac{1}{|\Aut_M(\Gamma)|}
\; A_{k(\Gamma),f(\Gamma)} ^{nor}\;
\prod_j t_j ^{v_j(\Gamma)}\;,
\end{split}
\end{equation}
where
$$
A_{k,f} ^{nor} = \sum_{i_1,\dots,i_k;\;
h_1,\dots,h_{f-1} = 1} ^N
\left\langle 
e_{i_1}^2 \cdots e_{i_k} ^2\cdot
e_{h_1}^* e_{h_1} \cdots e_{h_{f-1}} ^* e_{h_{f-1}}
\right\rangle\;,
$$
$\Aut_M(\Gamma)$ is the automorphism group of a
graph $\Gamma$ as a M\"obius graph, and $k(\Gamma)$ is
the cross-cap genus of a non-orientable 
surface $S_\Gamma$.
 We notice the sharp contrast
between $A_{g,f} ^{or}$ and $A_{k,f} ^{nor}$, which reflects
a particular choice of a presentation of the fundamental
group $\pi_1(S_\Gamma)$ of  a closed surface $S_\Gamma$.
Every  simple finite-dimensional real 
von Neumann algebra is a full matrix algebra
over either the reals $\mathbb{R}$ or quaternions $\mathbb{H}$.  
We recover the graphical expansion 
formulas for real symmetric and quaternionic self-adjoint
matrix integrals of \cite{BIPZ, GHJ, MW} from 
 (\ref{eq:introexpreal}). 
An explicit computation is also carried out for 
real Clifford algebras  \cite{Yu}.

Here we emphasize again that even though their expressions
look dependent
on a presentation of $\pi_1(S_\Gamma)$, the
quantity $A_{g,f} ^{or}$ is an invariant of an orientable surface
of topological type $(g,f)$, and $A_{k,f} ^{nor}$   is an invariant
of a non-orientable surface of 
topological type $(k,f)$.
When the von Neumann algebra $A$
in our theory is  simple, the
invariants $A_{g,f} ^{or}$ and $A_{k,f} ^{nor}$   
do not show any significance. 
The invariants become more interesting when the
algebra is complicated.
Now we notice that every finite-dimensional von Neumann algebra
is semi-simple, and hence is decomposable into simple factors. 
When we apply the decomposition of $A$ into simple factors 
 in the integral of (\ref{eq:introexp}) or
(\ref{eq:introexpreal}), due to the logarithm in front
of the integral, it becomes
the sum of the integral for each simple factor. 
Therefore, any topological invariant given as 
$A_{g,f} ^{or}$ or $A_{k,f} ^{nor}$ is computable in terms
of simple ones. 

This idea can be concretely carried out for the real or
complex group algebra of a finite group $G$. Using the complex
group algebra $\mathbb{C}[G]$, we obtain
\begin{multline}
\label{eq:introCGexp}
\log \int_{\{x\in \mathbb{C}[G]\;|\;x=x^*\}}
\exp\left(-\frac{1}{2}\;\rchi_\reg(  x^2 )\right)
\exp\left(\sum_{j} \frac{t_j}{j}\rchi_\reg(x^j )
\right) d\mu(x)\\
= \sum_{\substack{\Gamma  \text{ connected}\\
\text{ribbon graph}}}
\frac{1}{|\Aut_R \Gamma|}|G|^{\rchi(S_\Gamma)-1}
|\Hom(\pi_1(S_\Gamma),G)|\; 
\prod_{j} t_j ^{v_j(\Gamma)}\;,
\end{multline}
where $\rchi(S_\Gamma)$ is the Euler characteristic of $S_\Gamma$,
and $\rchi_\reg$ denotes the character of the regular representation 
of $G$ on $\mathbb{C}[G]$ linearly extended 
to the whole algebra.
Notice that the formula gives a generating function for the cardinality
of the representation variety $\Hom(\pi_1(S),G)$ of a closed oriented
surface $S$ in the group $G$. 
With the real group algebra $\mathbb{R}[G]$ of $G$, we have
\begin{multline}
\label{eq:introRGexp}
\log \int_{\{x\in \mathbb{R}[G]\;|\;x=x^*\}}
\exp\left(-\frac{1}{4}\;\rchi_\reg(  x^2 )\right)
\exp\left(\sum_{j} \frac{t_j}{2j}\rchi_\reg(x^j )
\right) d\mu(x)\\
= \sum_{\substack{\Gamma  \text{ connected}\\
\text{M\"obius graph}}}
\frac{1}{|\Aut_M \Gamma|}|G|^{\rchi(S_\Gamma)-1}
|\Hom(\pi_1(S_\Gamma),G)|\; 
\prod_{j} t_j ^{v_j(\Gamma)}\;.
\end{multline}
Surprisingly, the 
RHS of (\ref{eq:introRGexp}) has the same expression as in 
(\ref{eq:introCGexp}), with the only difference being replacing
ribbon graphs with M\"obius graphs. 
These generating functions were reported  in 
an earlier paper \cite{MY}.

Let $G$ be a finite group and $\Ghat$ the set of equivalence
classes of complex irreducible representations of $G$. 
The most fundamental formula in representation theory of
finite groups is the one that expresses the order of the group in terms
of a square sum of the dimensions of irreducible
representations of $G$:
\begin{equation}
\label{eq:classical}
|G| = \sum_{\lam\in\Ghat} (\dim \lam)^2 \;.
\end{equation}
The formula follows from the
decomposition of the group algebra
into irreducible fctors:
\begin{equation}
\label{eq:GA}
\mathbb{C}[G] = \bigoplus_{\lam\in\Ghat} \;\End(\lam)\;.
\end{equation}
In 1978, Mednikh \cite{Med}
discovered a remarkable generalization of the
classical formula  (\ref{eq:classical}):
\begin{equation}
\label{eq:Med}
\sum_{\lam\in\Ghat} (\dim \lam)^{\rchi(S)}
= |G|^{\rchi(S)-1}|\Hom(\pi_1(S), G)|\;,
\end{equation}
where $S$ is a  compact Riemann surface. 
When $S=S^2$, (\ref{eq:Med}) reduces to 
(\ref{eq:classical}).
Note that (\ref{eq:GA}) is a von Neumann algebra isomorphism.
Thus the integral of (\ref{eq:introCGexp}) over the self-adjoint
elements of $\mathbb{C}[G]$ becomes the  sum of hermitian
matrix integrals. 
It is now easy to see that evaluation of the integral of
(\ref{eq:introCGexp}) using (\ref{eq:GA})
yields Mednykh's formula (\ref{eq:Med}). 

For a non-orientable surface $S$, the formula for the number of
representations of $\pi_1(S)$ 
involves more detailed information on irreducible representations
of $G$. 
Using  the Frobenius-Schur
indicator of irreducible characters \cite{FS}, 
we decompose the set of complex irreducible representations $\Ghat$
into the union of three disjoint subsets, corresponding to
real, complex, and quaternionic irreducible representations:
\begin{equation}
\label{eq:FSindicator}
\begin{split}
\Ghat _1 &= \bigg\{\lam\in \Ghat \;\bigg|\; \frac{1}{|G|}\sum_{w\in G}
\rchi_\lam(w^2)  = 1\bigg\}\; ;\\
\Ghat _2 &= \bigg\{\lam\in \Ghat \;\bigg|\; \frac{1}{|G|}\sum_{w\in G}
\rchi_\lam(w^2)  = 0\bigg\}\; ;\\
\Ghat _4 &= \bigg\{\lam\in \Ghat \;\bigg|\; \frac{1}{|G|}\sum_{w\in G}
\rchi_\lam(w^2)  = -1\bigg\}\; .
\end{split}
\end{equation}
The suffix $1,2$ or $4$ indicates
the dimension of the base field $\mathbb{R}$, $\mathbb{C}$
or $\mathbb{H}$, respectively. 
In the fundamental paper of Frobenius and Schur \cite{FS}
published in 1906, 
we find
\begin{equation}
\label{eq:FScount}
\sum_{\lam\in\Ghat_1}
(\dim\lam)^{\rchi(S)} 
+
\sum_{\lam\in\Ghat_4}
(-\dim\lam)^{\rchi(S)}
= |G|^{\rchi(S)-1}|\Hom(\pi_1(S),G)|   \;.
\end{equation}
 It is somewhat strange that
a formula for non-orientable surfaces was known  
much earlier than its orientable counterpart. Actually, 
Frobenius and Schur obtained the formula as a counting formula
for the number of group elements satisfying
$x_1 ^2 \cdots x_k ^2 = 1$, but no relation to surface topology
was in their motivation.
 If we take $S=\mathbb{R}P^2$, then 
the formula reduces to the well-known formula
\cite{Isaacs, Serre}
$$
\sum_{\lam\in\Ghat_1} \dim\lam
- \sum_{\lam\in\Ghat_4} \dim\lam
= \text{ the number of involutions of  } G\;.
$$
The formula (\ref{eq:FScount}) immediately follows from
the generating function (\ref{eq:introRGexp}) and the decomposition
of $\mathbb{R}[G]$ into simple factors, which include real, complex,
and quaternionic matrix algebras.

In a beautiful paper of Pierre van Moerbeke \cite{VM}, we see the list
of matrix-type integrals and the nonlinear integrable systems that 
characterize the integrals as functions on the potential. The simple
von Neumann algebra integrals all fit into his scope. More general
von Neumann algebra integrals of this article can be 
considered as a multi-matrix model with trivial interaction terms
between matrices. They can be also interpreted as a matrix integral
over an algebra different from $\mathbb{R}, \mathbb{C}$ and
$\mathbb{H}$. In either point of view, we do not have any clear
picture on the relation between our formulas
(\ref{eq:introexp}), (\ref{eq:introexpreal}) and integrable equations. 
Since the generating functions for the Hurwitz numbers and the 
Gromov-Witten invariants of the Riemann surfaces are proven to
satisfy integrable systems \cite{OP1, OP2}, the integrability of the
von Neumann algebra integrals seems to pose a condition on the 
structure of the algebra. However, the present article does not address this
question.

The study of the volume of the representation variety 
$\Hom(\pi_1(S),G)$ of a surface group in a compact
connected simply connected semi-simple Lie group 
is carried out by many authors including
Witten \cite{W1991}, Gross-Taylor \cite{GrossTaylor},
 and Liu \cite{L1,L2,L3}. 
Although their focus was on the moduli space
$\mathcal{M}(S,G)$
of flat $G$-connections on a closed surface $S$, through a
relation 
$$
\mathcal{M}(S,G) = \frac{\Hom(\pi_1(S),G)}{G/Z(G)}\;,
$$
the study of moduli spaces is equivalent to that of representation 
varieties. Here $Z(G)$ is the center of the group $G$ that
acts trivially on the representation variety through conjugation.
It is interesting to note  that
exactly the same
formulas (\ref{eq:Med}) and (\ref{eq:FScount}) 
hold for a compact Lie group if the infinite sum
of LHS converges absolutely and the cardinality is interpreted
as the \emph{volume} of the variety in an appropriate sense.  
The naive extension of the method of this article does not 
work for the case of Lie groups, however, because the
von Neumann algebra involved becomes infinite-dimensional and
the integration (\ref{eq:introCGexp}) makes no sense.

The present paper is organized as follows. 
We review the notion of ribbon graphs and \emph{M\"obius}
graphs in Section~\ref{sect:graph}. Then
in Section~\ref{sect:vna}, we compute 
 Feynman diagram expansion of integrals
over a finite-dimensional complex von Neumann algebra
in terms
of ribbon graphs.
Since integrals over a real von Neumann algebra behave
differently, they are treated separately in Section~\ref{sect:realvna}.
The generating functions for the number of representations
of surface groups in a finite group are given in 
Section~\ref{sect:generating}. As an application, 
 we give a new proof of 
the formulas of Frobenius-Schur and Mednykh.

\begin{ack}
The authors thank
Michael Penkava and Andrew Waldron for
discussions on non-commutative matrix integrals.
\end{ack}

\bigskip

\section{Ribbon graphs and M\"obius graphs}
\label{sect:graph}

In this section, we review basic facts about graphs
drawn on an orientable or non-orientable surface.
Definition of the automorphism groups of these
graphs is crutial
when we use them to compute 
non-commutative matrix integrals. 
Graphs on an oriented surface are called \emph{ribbon}
or \emph{fat} graphs. We refer to 
\cite{BIZ, Harer1986, Harer-Zagier, Kontsevich,   O, 
OP1, OP2, Penner, W1991a}
for the use of ribbon/fat graphs in the study of
moduli spaces of Riemann surfaces and related topics.
The double line notation was first introduced 
in \cite{'tHooft}, generalizing the graphical expansion
idea of \cite{Feynman}. Graphs on 
complex algebraic curves
were studied from the quite different point of view of 
Grothendieck's \emph{dessins d'enfants} 
\cite{Belyi, Grothendieck, Schneps, SL}. A relation 
between Strebel differentials \cite{Gardiner, Strebel}
and dessins d'enfants was studied in \cite{MP1998}.
The terminology \emph{M\"obius graph} was 
introduced in \cite{MW}
for a graph on a surface that is not oriented in order
to avoid possible confusion, since ribbon or fat graphs 
are usually oriented. Graphs on an orientable or
non-orientable surface are also called \emph{maps}.
Maps were studied mainly in the context of 
map coloring theorems \cite{GT, Ringel}.

\begin{figure}[htb]
\centerline{\epsfig{file=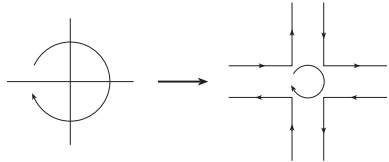, width=2in}}
\caption{A vertex with a cyclic order given to incident
half-edges. It is placed on a plane with the clockwise 
orientation. The half-edges become crossroads with a
compatible orientation at the boundary.}
\label{fig:orvertex}
\end{figure}

For a detailed treatment of ribbon graphs, we refer to the
references cited above. Only a brief description is given here.
A graph $\Gamma = (\mathcal{V}, \mathcal{E}, \iota)$
consists of a set of vertices $\mathcal{V}$, a set of 
edges $\mathcal{E}$, and an incidence relation
$$
\iota: \mathcal{E}\longrightarrow (\mathcal{V}\times
\mathcal{V})\big/ \mathfrak{S}_2.
$$
Following \cite{MP1998}, let us introduce the \emph{edge
refinement} $\Gamma_\mathcal{E}$
of a graph $\Gamma$, which is the original 
graph together with a two-valent vertex 
(the midpoint) chosen from each edge.  A \emph{half-edge} of 
$\Gamma$ is an edge of its edge refinement. 
A \emph{ribbon graph} is a graph with a cyclic
order assigned at each vertex
 to the set of half-edges incident to the vertex.
When a cyclic order is given, 
a vertex can be placed on an oriented plane, and half-edges
incident to the vertex can be represented by double lines.
The orientation of the plane gives an orientation of the
ribbon-like structure, and its boundaries inherit a compatible
orientation (see Figure~\ref{fig:orvertex}).

\begin{figure}[htb]
\centerline{\epsfig{file=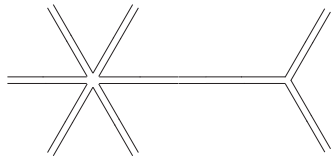, width=1.7in}}
\caption{A ribbon graph is obtained by connecting
cyclically ordered vertices with a ribbon like edge
preserving the orientation.}
\end{figure}

A topological realization of a
ribbon graph $\Gamma$ is obtained by connecting these half-edges
in an orientation-compatible manner. Since each 
boundary has a well-defined orientation, we can attach
an oriented disk to the boundary and form a compact oriented
surface $S_\Gamma$. Let $f(\Gamma)$ denote the number
of disks attached. This number is uniquely determined by the
ribbon graph structure of a graph. 
The attached disks, together with the vertices and edges of
$\Gamma$, form a \emph{cell-decomposition} of the surface
$S_\Gamma$. The genus of the surface is determined by 
the formula for the Euler characteristic
\begin{equation}
\label{eq:Euler}
\rchi(S_\Gamma)=2-2g(S_\Gamma)=v(\Gamma)-e(\Gamma)
+f(\Gamma),
\end{equation}
where $v(\Gamma) = |\mathcal{V}|$ is the number of
vertices and $e(\Gamma)=|\mathcal{E}|$ the number of
edges of $\Gamma$. 

\begin{figure}[htb]
\centerline{\epsfig{file=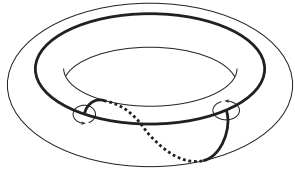, width=1.5in}}
\caption{A graph drawn on an oriented surface. At each
vertex, the orientation of the surface determines a
cyclic order of the edges incident to the vertex.}
\label{fig:surfacegraph}
\end{figure}

Conversely, if a connected
graph $\Gamma$ is drawn on an oriented surface
$S$ in a way that $S\setminus \Gamma$ is the union of
disjoint open disks, then $\Gamma$ is a ribbon graph that
defines a cell-decomposition
of $S$. The cyclic order of half-edges incident to a vertex
is determined by the orientation of the surface
(see Figure~\ref{fig:surfacegraph}).
Obviously we have $S=S_\Gamma$. 

\begin{Def}[\cite{MP1998}]
\label{def:ribbonauto}
Let $\Gamma$ be a ribbon graph.
The group $\Aut_R \Gamma$ of automorphisms of
$\Gamma$ consists of graph automorphisms of
the edge refinement $\Gamma_\mathcal{E}$ that 
preserves the cyclic order at each vertex of $\Gamma$.
\end{Def}

\begin{figure}[htb]
\centerline{\epsfig{file=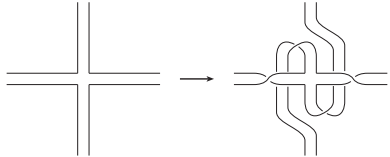, width=2in}}
\caption{A vertex flip operation. It reverses the cyclic
order at a vertex, and gives an extra twist to each half-edge incident 
to the vertex.}
\label{fig:vertexflip}
\end{figure}

In a ribbon graph, an edge connects two oriented vertices
 in the orientation-compatible manner. If we 
connect vertices without paying attention to  the orientation, 
then we obtain
a \emph{M\"obius} graph. An edge connecting two oriented
vertices is \emph{not twisted} if the connection is consistent 
with the orientation, and is \emph{twisted} otherwise. 
Thus a double twist is the same as no twist. A new operation 
allowed in a M\"obius graph that preserves the M\"obius
graph structure is a \emph{vertex flip} at a vertex. 
This operation reverses the cyclic order assigned at the
vertex, and twists every half-edge incident to the vertex
(see Figure~\ref{fig:vertexflip}). If an edge is incident to
a vertex and forms a loop, then the vertex flip at this vertex
does not change the twist of the edge.

We can formalize the definition of a M\"obius graph
in the following way.

\begin{Def}
\label{def:moebius}
A \emph{M\"obius graph} is the equivalence class of
  ribbon graphs with a $\mathbb{Z}/2\mathbb{Z}$-color
assigned to each edge. Two edge-colored ribbon graphs
are equivalent if one is obtained from the other by 
a sequence of vertex flip operations. A vertex
flip  reverses 
the cyclic order of a vertex and the color of
each half-edge incident to it. The group 
$\Aut_M \Gamma$ of automorphisms of a M\"obius graph
consists of graph automorphisms of the edge refinement
of the underlying graph $\Gamma$ that preserve the 
equivalence class of the edge-colored ribbon graph.
\end{Def}

\begin{figure}[htb]
\centerline{\epsfig{file=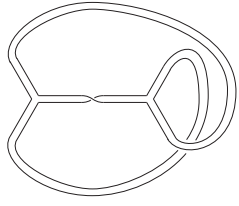, width=1.2in}}
\caption{A M\"obius graph.}
\end{figure}

A topological realization of a M\"obius graph 
is the realization of the $\mathbb{Z}/2\mathbb{Z}$-color
of each edge as a twist or non-twist. Each boundary component
of a M\"obius graph is a circle, without any 
consistent orientation. By attaching an open disk to each
boundary circle, a M\"obius graph gives rise to a closed surface
without orientation. Let us denote this surface by
$S_\Gamma$ and by $f(\Gamma)$ the number
of disks, as before. We note that $\Gamma$ defines
a cell-decomposition of $S_\Gamma$. 
Every closed non-orientable surface is constructed
by removing $k$ disks from a sphere and attaching 
a cross-cap at each hole. The number $k$ is the
\emph{cross-cap genus} of the surface, and the Euler
characteristic of the surface is given by $2-k$. 
If the surface $S_\Gamma$ is non-orientable, then
we have
\begin{equation}
\label{eq:nonoreuler}
\rchi(S_\Gamma) = 2-k(S_\Gamma) = v(\Gamma)-e(\Gamma)
+f(\Gamma).
\end{equation}

A ribbon graph $\Gamma$ is also
 a M\"obius graph. If $\Gamma$
and its flip $\Gamma^t$ (the graph obtained by applying the vertex flip
operation at every vertex simultaneously)
is the same ribbon graph, then we have
\begin{equation}
\label{eq:flipgraph2}
|\Aut_M\Gamma|=2\; |\Aut_R \Gamma|.
\end{equation}
Otherwise, $\Gamma$ and $\Gamma^t$ are
different ribbon graphs but the same as a M\"obius graph,
and we have
\begin{equation}
\label{eq:flipgraph}
\Aut_M\Gamma \isom  \Aut_R \Gamma.
\end{equation}

\bigskip

\section{Integrals over a finite-dimensional 
complex von Neumann
algebra}
\label{sect:vna}

In this section we define the integrals
over a finite-dimensional complex
von Neumann algebra
that we study, and establish their graphical 
expansion formulas in terms of ribbon graphs.

\begin{Def}
\label{def:vna}
A finite-dimensional  \emph{complex von Neumann
algebra} is a finite-dimensional
$\mathbb{C}$-algebra with a conjugate-linear anti-isomorphism
$*:A\rightarrow A$ and a $\mathbb{C}$-linear
map called \emph{trace} $\langle\; \rangle:A\rightarrow
\mathbb{C}$ that satisfy the following conditions
for every $a,b\in A$:
\begin{equation}
\begin{split}
\label{eq:vna}
(a^*)^* &= a\\
(ab)^* &= b^*a^*\\
\langle a^*\rangle &= \overline{\langle a\rangle}\\
\langle ab\rangle &= \langle ba\rangle\\
\langle 1\rangle &=1\\
\langle aa^* \rangle &>0\;, \qquad a\ne 0\;.
\end{split}
\end{equation}
If $A$ is an $\mathbb{R}$-algebra with a real valued 
trace, then it is called a \emph{real} von Neunamm algebra.
\end{Def}

To avoid confusion, we only deal with complex
von Neumann algebras in this section. Real ones
are considered in Section~\ref{sect:realvna}.
A finite-dimensional von Neumann algebra $A$ is a real
vector space with a non-degenerate 
hermitian inner product defined by
\begin{equation}
\label{eq:hermitianinner}
\langle a,b\rangle = \langle ab^*\rangle\;.
\end{equation}
As usual, an invertible linear transformation 
of $A$ that preserves
the hermitian form is called a unitary transformation. 
We denote by $U(A)$ the group of unitary transformations
of $A$. 
A real vector subspace of $A$ consisting of 
self-adjoint elements
\begin{equation}
\label{eq:self-adjoint}
\mathcal{H}_A = \{a\in A\;|\;a^* = a\}
\end{equation}
is of our particular interest. We note that $\mathcal{H}_A$
is invariant under the conjugation action of $U(A)$. 
Let us denote by $dx$ the translation invariant
Lebesgue measure of the real
vector space $\mathcal{H}_A$ that is also invariant under
the conjugation action of $U(A)$. 
We notice that the quadratic form $\langle x^2\rangle$
is positive definite on the space $\mathcal{H}_A$
of  self-adjoint elements. We denote by
\begin{equation}
\label{eq:dmux}
d\mu (x) = 
\frac{dx}{\int_{\mathcal{H}_A}
e^{-\frac{1}{2}\langle x^2\rangle }\; dx}
\end{equation}
 the normalized Lebesgue measure on $\mathcal{H}_A$.

Our subject of study
is the following integral 
\begin{equation}
\label{eq:HAintegral}
\prod_{j=1} ^n \frac{1}{v_j ! \cdot j^{v_j}}
\int_{\mathcal{H}_A} e^{-\frac{1}{2}\langle x^2\rangle}
\prod_{j=1} ^n \langle x^j\rangle^{v_j}d\mu(x)
\end{equation}
for every $n$-tuple of positive integers
$(v_1,\dots,v_n)\in\mathbb{N}^n$, $n=0,1,2,\cdots$. 
The constant factor in front of the integral is placed
for a combinatorial reason explained later in this section.
To consider a generating function
of these integrals, it is more convenient to introduce 
\begin{equation}
\label{eq:eofv}
\mathbf{e}(v_1,v_2,v_3,\cdots) =  \sum_{j\ge 1} j v_j
\end{equation}
 and the sum 
of the integrals over all elements of
$$
\mathbb{N}^\infty = 
\lim_{\substack{\longrightarrow\\ n}}
\mathbb{N}^n
=\{(v_1,v_2,v_3,\dots)\;|\;v_j=0 \;\text{ for }\;  j>>0\}
$$
with a fixed value of $\mathbf{e}(v_1,v_2,v_3,\cdots)$. 
Notice that for every finite value of $n$, 
$$
(v_1,v_2,v_3,\dots)\in\mathbb{N}^n
\qquad{\text{if}}\qquad 
\mathbf{e}(v_1,v_2,v_3,\cdots)\le n\;.
$$
Thus let us define
\begin{equation}
\label{eq:HAgenerating}
Z_A ^\mathbb{C} (t_1,t_2,t_3,\dots) = 
\sum_{n=0} ^\infty 
\sum_{\substack{
(v_1,v_2,v_3,\dots)\in\mathbb{N}^\infty\\
\mathbf{e}(v_1,v_2,v_3,\cdots)=n}}
\prod_{j\ge1} ^{\fin} \frac{t_j ^{v_j}}{v_j ! j^{v_j}}
\int_{\mathcal{H}_A} e^{-\frac{1}{2}\langle x^2\rangle}
\prod_{j\ge1} ^{\fin} \langle x^j\rangle^{v_j}d\mu(x)\;,
\end{equation}
where $t_1,t_2,t_3,\dots$ are expansion parameters
carrying the weight 
\begin{equation}
\label{eq:degtj}
\deg t_j = j\;.
\end{equation}
The monomial $\prod_j t_j ^{v_j}$ for every 
$(v_1,v_2,v_3,\cdots)$ satisfying
$\mathbf{e}(v_1,v_2,v_3,\cdots)=n$ has weighted homogeneous 
degree $n$ by (\ref{eq:eofv}) and (\ref{eq:degtj}). 
Hence (\ref{eq:HAgenerating}) is an infinite sum of 
weighted homogeneous polynomials of degree $n$ for every 
$n\ge 0$.

Symbolically, we can write the generating function
in an integral form
\begin{equation}
\label{eq:HAgenint}
Z_A  ^\mathbb{C}(t) = 
\int_{\mathcal{H}_A} e^{-\frac{1}{2}\langle x^2\rangle}
e^{\sum_{j=1} ^\infty
\frac{t_j}{j} \langle x^j\rangle^{v_j}}d\mu(x)\;.
\end{equation}
As an actual integral, (\ref{eq:HAgenint}) is ill
defined because of the infinite sum in the exponent.
There is a way to make it well-defined so that 
(\ref{eq:HAgenerating}) is a rigorous asymptotic
expansion of (\ref{eq:HAgenint}). Since we do not
employ this point of view in this paper, we refer to 
\cite{Mulase1995,Mulase1998}
for more detail, and work on the expansion 
form only. 

Let   $\{e_1,\dots,e_N\}$ be an orthonormal
basis for $A$ with respect to the hermitian form 
(\ref{eq:hermitianinner}), where
$N=\dim_{\mathbb{C}} A$. Since
$$
\langle e_i,e_j\rangle = \langle e_i e_j ^*\rangle
= \langle e_j ^* (e_i ^*)^*\rangle
= \langle e_j ^*,e_i ^*\rangle\;,
$$
$\{e_1 ^*,\dots,e_N ^*\}$ also forms an orthonormal
basis for $A$. For every $a\in A$ we have
\begin{equation}
\label{eq:eiej}
a = \sum_{j=1} ^N \langle a,e_j \rangle \;e_j 
=\sum_{j=1} ^N \langle a,e_j ^*\rangle \;e_j ^*\;.
\end{equation}
Equivalently, 
\begin{equation}
\label{eq:ab}
\langle a,b\rangle = 
\sum_{j=1} ^N \langle a,e_j \rangle \;\langle e_j ,b\rangle
=\sum_{j=1} ^N \langle a,e_j ^*\rangle \;\langle e_j ^*, b\rangle
\end{equation}
holds for every $a,b\in A$.

\begin{lem}
\label{lem:exy}
Choose two elements 
\begin{equation*}
x = \sum_{i=1} ^N x_i \;e_i 
\qquad\text{and}\qquad
y = \sum_{i=1} ^N y_i \;e_i
\end{equation*}
of $A$,
and consider $e^{\langle xy\rangle}$
as a function in $2N$ variables 
$$
(x_1,\dots,x_N,y_1,\dots,y_N)\in \mathbb{C}^{2N}\;.
$$
With respect to the differential operator 
\begin{equation}
\label{eq:ddy}
\frac{\partial}{\partial y} = \sum_{i=1} ^N 
\frac{\partial}{\partial y_i}\;e_i ^*\;,
\end{equation}
we have
\begin{equation}
\label{eq:ddyexy}
\frac{\partial}{\partial y}\; e^{\langle xy\rangle}
= x\; e^{\langle xy\rangle}\;.
\end{equation}
In particular,
\begin{equation}
\label{eq:ddytrace}
\left\langle \left( \frac{\partial}{\partial y}
\right) ^j\right\rangle ^m\; e^{\langle xy\rangle}
= \langle x^j\rangle ^m e^{\langle xy\rangle}
\end{equation}
for every $j,m>0$. 
\end{lem}

\begin{proof}
By definition, 
\begin{equation*}
\begin{split}
\frac{\partial}{\partial y}\; e^{\langle xy\rangle}
&= \sum_i \frac{\partial}{\partial y_i}\; e_i ^*
\exp\left( \left\langle \sum_j x_j\; e_j \sum_k y_k \; e_k
\right\rangle\right)\\
&=\sum_i \sum_j x_j \; \langle e_j e_i\rangle \; e_i ^*\;
e^{\langle xy\rangle}\\
&=\sum_j x_j \sum_i\langle e_j,e_i ^* \rangle \; e_i ^*\;
e^{\langle xy\rangle}\\
&= \sum_j x_j\; e_j\; e^{\langle xy\rangle}\\
&=x\; e^{\langle xy\rangle}\;.
\end{split}
\end{equation*}
Using the linearity of
the trace and (\ref{eq:ddyexy}) repeatedly, we obtain
(\ref{eq:ddytrace}).
\end{proof}

\begin{lem}
\label{lem:complexlaplace}
Let $A$ be a finite-dimensional complex von Neumann algebra.
Then we have the following Laplace transform formula for
(\ref{eq:HAintegral}):
\begin{equation}
\label{eq:complexlaplace}
\int_{\mathcal{H}_A} e^{-\frac{1}{2}\langle x^2\rangle}
\prod_{j=1} ^n \langle x^j\rangle^{v_j}d\mu(x)
=\left. \prod_{j=1} ^n \left\langle \left(
\frac{\partial}{\partial y}\right)^j\right\rangle^{v_j}
e^{\frac{1}{2}\langle(y+y^*)^2\rangle}\right|_{y=0}\;.
\end{equation}
\end{lem}

\begin{proof}
For $y\in A$ of Lemma~\ref{lem:exy}, its adjoint is given by
$$
y^* = \sum_{i=1} ^N \overline{y_i}\; e_i ^*\;.
$$
Note that $\partial y^*/\partial y=0$
since $\partial \overline{y_i}/\partial y_j = 0$ for any
$i$ and $j$. 
Now Lemma~\ref{lem:exy} yields 
\begin{equation*}
\int_{\mathcal{H}_A} e^{-\frac{1}{2}\langle x^2\rangle}
 \langle x^j\rangle^{m}d\mu(x)
= \left.\left\langle \left(
\frac{\partial}{\partial y}\right)^j\right\rangle^{m}
\int_{\mathcal{H}_A} e^{-\frac{1}{2}\langle x^2\rangle}
e^{\langle x(y+y^*)\rangle} d\mu(x)\right|_{y=0}\;.
\end{equation*}
Since $y+y^*\in \mathcal{H}_A$ and $d\mu(x)$ is a
translational invariant  measure, we have
\begin{equation*}
\begin{split}
\int_{\mathcal{H}_A} e^{-\frac{1}{2}\langle x^2\rangle}
e^{\langle x(y+y^*)\rangle} d\mu(x)
&=
\int_{\mathcal{H}_A} e^{-\frac{1}{2}\langle (x-(y+y^*))^2\rangle}
e^{\frac{1}{2}\langle (y+y^*)^2\rangle} d\mu(x)\\
&=e^{\frac{1}{2}\langle (y+y^*)^2\rangle}\;.
\end{split}
\end{equation*}
Eqn.~(\ref{eq:complexlaplace}) follows from these formulas. 
\end{proof}

\begin{figure}[htb]
\centerline{\epsfig{file=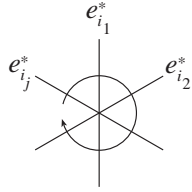, width=1in}}
\caption{A $j$-valent vertex with a cyclic order given to incident
half-edges.
Each half-edge is labeled by $i_k$, and an element
$e_{i_k} ^*$ is assigned.}
\label{fig:jpower}
\end{figure}

In the same way as in the proof of Lemma~\ref{lem:exy},
we obtain
\begin{equation*}
\begin{split}
\frac{\partial}{\partial y} \;e^{\frac{1}{2}\langle
(y+y^*)^2\rangle} &= 
\sum_i \frac{\partial}{\partial y_i}\; e_i ^*\;
e^{\frac{1}{2}\langle (\sum_j y_j\; e_j + \sum_j \overline{y_j}\; e_j ^*)^2
\rangle}\\
&=\sum_i \langle (y+y^*) \;e_i \rangle \; e_i ^*\;
e^{\frac{1}{2}\langle(y+y^*)^2\rangle}\\
&=\sum_i \langle (y+y^*) ,e_i ^* \rangle \; e_i ^*\;
e^{\frac{1}{2}\langle(y+y^*)^2\rangle}\\
&=(y+y^*) \; e^{\frac{1}{2}\langle(y+y^*)^2\rangle}\;.
\end{split}
\end{equation*}
In particular, 
\begin{equation*}
\frac{\partial}{\partial y_i}
e^{\frac{1}{2}\langle(y+y^*)^2\rangle}
= \langle (y+y^*) ,e_i ^* \rangle \;
e^{\frac{1}{2}\langle(y+y^*)^2\rangle}\;,
\end{equation*}
and hence
\begin{equation}
\label{eq:dyidyjcomplex}
\begin{split}
\left.\frac{\partial}{\partial y_i}\frac{\partial}{\partial y_j}
e^{\frac{1}{2}\langle(y+y^*)^2\rangle}\right|_{y=0}
&= 
\left.\frac{\partial}{\partial y_i}\langle (y+y^*) ,e_j ^* \rangle 
\right|_{y=0}\\
&=\langle e_i ,e_j ^*\rangle\\
&=\langle e_i e_j \rangle\;.
\end{split}
\end{equation}
Our purpose is to compute
\begin{equation}
\label{eq:complexdiff}
 \sum_{\substack{
(v_1,v_2,v_3,\dots)\in\mathbb{N}^\infty\\
\mathbf{e}(v_1,v_2,v_3,\cdots)=n}}
\prod_{j\ge1} ^{\fin} 
\frac{1}{v_j !\cdot j^{v_j}}
\left.\left\langle \left(
\frac{\partial}{\partial y}\right)^j\right\rangle^{v_j}
e^{\frac{1}{2}\langle(y+y^*)^2\rangle}\right|_{y=0}\;.
\end{equation}
To this end, first we observe
\begin{equation}
\label{eq:jpower}
\left\langle\left(\frac{\partial}{\partial y}\right)^j\right\rangle
=\sum_{i_1,\dots,i_j } 
\frac{\partial}{\partial y_{i_1}}\cdots
\frac{\partial}{\partial y_{i_j}}
\big\langle e_{i_1} ^*\cdots e_{i_j} ^*\big\rangle\;.
\end{equation}
Notice  that 
$\big\langle e_{i_1} ^*\cdots e_{i_j} ^*\big\rangle$
is invariant under cyclic permutations. 
For every factor (\ref{eq:jpower}) of 
(\ref{eq:complexdiff}), let us assign a $j$-valent vertex
with $j$ half-edges incident to it, with a cyclic order of these
half-edges. Every half-edge corresponds to an index $i_k$,
$k=1,\cdots,j$, and we assign $e_{i_k} ^*\in A$ to this
half-edge. We then assign 
$\big\langle e_{i_1} ^*\cdots e_{i_j} ^*\big\rangle$
to this vertex (see Figure~\ref{fig:jpower}).

\begin{figure}[htb]
\centerline{\epsfig{file=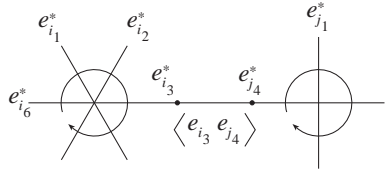, width=2in}}
\caption{A half-edge labeled by $i_3$ of the left vertex
is connected with a half-edge labeled by $j_4$ of the
vertex at the right. A propagator
$\langle e_{i_3} e_{j_4}\rangle$ is assigned to this edge.}
\label{fig:edgeprop}
\end{figure}

For every $j=1,2,3,\cdots$, we draw $v_j$ $j$-valent vertices
with $\big\langle e_{i_1} ^*\cdots e_{i_j} ^*\big\rangle$
assigned. Every vertex has $j$ degrees of freedom coming
from cyclic rotations. This redundancy is compensated
by the factor $j^{v_j}$ in (\ref{eq:complexdiff}). 
The redundancy of permuting $v_j$ vertices of the same
valence $j$ is compensated by $v_j !$ in (\ref{eq:complexdiff}). 
To indicate the effect of (\ref{eq:dyidyjcomplex}),
we connect two half-edges according to the paired
differentiation. We notice that since we set $y=0$ after
differentiation, no term in (\ref{eq:complexdiff})
survives unless all differentiations are paired as in
(\ref{eq:dyidyjcomplex}). 
When we connect a half-edge labeled by $i_k$ at 
one vertex with another half-edge labeled by $h_\ell$, 
we assign $\langle e_{i_k}e_{h_\ell}\rangle$ to this
edge (see Figure~\ref{fig:edgeprop}).
This quantity is called the \emph{propagator}
of the edge. Notice that the propagator is symmetric
\begin{equation}
\label{eq:sympropagator}
\langle e_{i_k}e_{h_\ell}\rangle
= \langle e_{h_\ell}e_{i_k}\rangle\;,
\end{equation}
and hence we do not have any particular \emph{direction}
on our edge.

Here we have to be cautious when
two vertices are connected. For example, the connection
described in Figure~\ref{fig:edgeprop} preserves
the cyclic orders of two vertices. An edge is connecting
two vertices in the  \emph{orientation-preserving}
manner if the cyclic orders of the two vertices agree
when the edge is shrunk to a point and the two vertices are
put together. Otherwise, the edge is 
\emph{orientation-reversing}. 
All connections we make in this section 
should be orientation-preserving.
When all half-edges are paired and connected
in the orientation-preserving manner, we obtain
a \emph{ribbon graph} $\Gamma$. 
It is easy to see that the compensation of 
rotations around each vertex and permutations
of vertices of the same valence leads to the factor
of $1/|\Aut_R(\Gamma)|$ coming from the
automorphism of $\Gamma$
\cite{Mulase1998, MP1998}. 
The quantity $v_j$ represents the number of $j$-valent
vertices of $\Gamma$ by the construction. 
Thus
\begin{equation}
\label{eq:vandefromvj}
v(\Gamma) = \sum_j v_j\qquad
\text{and}
\qquad 
e(\Gamma) = \frac{1}{2}\; \sum_j jv_j
\end{equation}
represent the total number of vertices and edges of
$\Gamma$, respectively. Notice the combinatorial
constraint
$$
2e(\Gamma) = \mathbf{e}(v_1,v_2,v_3,\dots)\;,
$$
which comes from the fact that unless every half-edge is
paired with another one to
form an edge of a ribbon graph, the corresponding contribution
of $(v_1,v_2,v_3,\dots)$
in the sum of (\ref{eq:complexdiff}) is $0$. 
Summarizing, we have

\begin{prop}
\label{prop:ribbonexp}
Let $RG(e)$ denote the set of all ribbon
graphs, may or may not be connected,
 consisting of a total of $e$ edges. The number
of $j$-valent vertices of $\Gamma\in RG(e)$ is
denoted by $v_j = v_j(\Gamma)$.
For each $j$-valent vertex of   $\Gamma$, let
us assign 
$\big\langle e_{i_1} ^*\cdots e_{i_j} ^*\big\rangle$.
For every edge we assign $\langle e_{i_k}e_{h_\ell}\rangle$,
so that the incidence relation is consistent with the
relation described above, namely, this edge
connects the half-edge labeled by $i_k$ of a vertex
to the half-edge labeled by $h_\ell$ from another
vertex, which could be the same vertex.
Let $A_\Gamma ^{or}$ denote the sum 
with respect to  all indices of the product
of all contributions from vertices and edges. Then we have
\begin{equation}
\label{eq:graphcontribution}
\begin{split}
& \quad\sum_{\substack{
(v_1,v_2,v_3,\dots)\in\mathbb{N}^\infty\\
\mathbf{e}(v_1,v_2,v_3,\cdots)=2e}}
\int_{\mathcal{H}_A} e^{-\frac{1}{2}\langle x^2\rangle}
\prod_{j\ge1} ^{\fin}  \frac{1}{v_j !\cdot j^{v_j}}
 \langle x^j\rangle^{v_j}d\mu(x)\\
&=  \sum_{\substack{
(v_1,v_2,v_3,\dots)\in\mathbb{N}^\infty\\
\mathbf{e}(v_1,v_2,v_3,\cdots)=2e}}
\prod_{j\ge1} ^{\fin}  \frac{1}{v_j !\cdot j^{v_j}}
\left.\left\langle \left(
\frac{\partial}{\partial y}\right)^j\right\rangle^{v_j}
e^{\frac{1}{2}\langle(y+y^*)^2\rangle}\right|_{y=0}\\
&= \sum_{\Gamma\in RG(e)}
\frac{1}{|\Aut_R (\Gamma)|} \;A_\Gamma ^{or}\;,
\end{split}
\end{equation}
where $\Aut_R (\Gamma)$ is
 the automorphism group of ribbon graph $\Gamma$
defined in Section~\ref{sect:graph}. 
\end{prop}

Therefore, to evaluate the integral, it suffices to calculate
$A_\Gamma ^{or}$ for each ribbon graph $\Gamma$. A key 
fact is the following.

\begin{lem}
\label{lem:topology}
Let $\Gamma$ be a connected ribbon graph with two or
more vertices, and $E$ an edge of $\Gamma$
incident to two distinct vertices. Then the contribution
of the graph $A_\Gamma ^{or}$ is invariant under the edge-contraction:
\begin{equation}
\label{eq:topology}
A_\Gamma ^{or}= A_{\Gamma/ E} ^{or}\;,
\end{equation}
where $\Gamma/ E$ denotes the ribbon graph obtained
by shrinking $E$ to a point in $\Gamma$ and joining the
two incident vertices together.
\end{lem}

\begin{rem}
This invariance is  found in many
literatures including
\cite{W1991}. Witten uses the invariance to calculate 
quantum Yang-Mills theory over a Riemann surface
by approximation through 
lattice gauge theory. It appears also
in \cite{MW}. 
\end{rem}

\begin{proof}
Let $V_1$ and $V_2$ be the two vertices of $\Gamma$ 
incident to $E$. The contribution from $V_1$ 
can be written as $\langle a e_i ^*\rangle$ and 
that from $V_2$ as $\langle e_j ^* b\rangle$, where
$a$ and $b$ are products of the basis elements $e_k ^*$ of
the von Neumann algebra $A$. The invariance of the
edge contraction is local, and comes down to the
following computation:
\begin{equation}
\label{eq:edgecontraction}
\begin{split}
\sum_{i,j} \langle a e_i ^*\rangle \langle e_i e_j\rangle
\langle e_j ^* b\rangle 
&=
\sum_{i,j} \langle a , e_i \rangle \langle e_i, e_j ^*\rangle
\langle e_j ^* ,b^*\rangle\\
&=\langle a ,b^*\rangle\\
&=\langle a b\rangle\;.
\end{split}
\end{equation}
The quantity $\langle a b\rangle$ is exactly the
contribution of the new vertex obtained by joining
$V_1$ and $V_2$. 
\end{proof}

\begin{figure}[htb]
\centerline{\epsfig{file=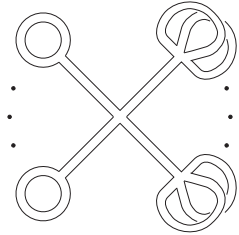, width=1.2in}}
\caption{A standard graph of topological type
$(g,f)$. It has $f-1$ tadpoles on the left, and $g$
bi-petal flowers on the right.}
\label{fig:standard}
\end{figure}

Every connected ribbon graph $\Gamma$ gives rise 
to an oriented surface $S_\Gamma$ whose Euler characteristic 
is determined by 
$$
\rchi(S_\Gamma) = 2-2g(S_\Gamma)
= v(\Gamma) -e(\Gamma) +f(\Gamma)\;.
$$
The graph defines  a cell-decomposition of 
$S_\Gamma$. The topological type of 
$\Gamma$ is $(g,f)$, the genus of the surface and the
number of faces of its cell-decomposition. 
Note that since the edge contraction operation decreases
$v(\Gamma)$ and $e(\Gamma)$ by one and preserves
the number of faces, the topological type is
preserved. A theorem of topology states that if
$\Gamma_1$ and $\Gamma_2$ are two connected
ribbon graphs with the same topological type,
 then by consecutive
applications of edge contraction and its inverse operation
(edge expansion), $\Gamma_1$ can be brought to 
$\Gamma_2$ \cite{Hatcher}. Therefore, to compute
$A_\Gamma ^{or}$, we can use our favorite graph of the
same topological type, for example,
a graph of Figure~\ref{fig:standard}.

\begin{prop}
\label{prop:standardcont}
Let $A_{g,f} ^{or}$ denote the contribution of the 
standard graph of Figure~\ref{fig:standard}. Then 
\begin{equation}
\label{eq:Agf}
A_{g,f} ^{or} = 
\sum_{\substack{i_1,\dots,i_g;j_1,\dots,j_g\\
h_1,\dots,h_{f-1} }} ^N
\left\langle
e_{i_1}e_{j_1}e_{i_1} ^* e_{j_1} ^*
\cdots e_{i_g}e_{j_g}e_{i_g} ^* e_{j_g} ^*\cdot 
e_{h_1}e_{h_1} ^*\cdots e_{h_{f-1}} e_{h_{f-1}} ^*
\right\rangle\;.
\end{equation}
\end{prop}

\begin{proof}
By definition, 
\begin{equation*}
\begin{split}
A_{g,f} ^{or} = 
\sum_{\substack{i_1,\dots,i_g;j_1,\dots,j_g\\
a_1,\dots,a_g;b_1,\dots,b_g\\
k_1,\dots,k_{f-1}; h_1,\dots,h_{f-1} }} ^N
&\left\langle
e_{a_1} ^*e_{b_1} ^* e_{i_1} ^* e_{j_1} ^*
\cdots e_{a_g} ^*e_{b_g} ^*e_{i_g} ^* e_{j_g} ^*\cdot 
e_{h_1} ^* e_{k_1} ^*\cdots e_{h_{f-1}}^* e_{k_{f-1}} ^*
\right\rangle\\
&\times
\langle e_{a_1} e_{i_1}\rangle \; \langle e_{b_1} e_{j_1}\rangle
\cdots
\langle e_{a_g} e_{i_g}\rangle \; \langle e_{b_g} e_{j_g}\rangle
\cdot
\langle e_{h_1} e_{k_1}\rangle \; \cdots 
\langle e_{h_{f-1}} e_{k_{f-1}}\rangle\;.
\end{split}
\end{equation*}
Using cyclic invariance of the trace and (\ref{eq:ab}), 
the desired formula (\ref{eq:Agf}) follows.
\end{proof}

The generating function $Z_A ^\mathbb{C}(t)$ of 
(\ref{eq:HAgenerating}) is expanded in terms of all
ribbon graphs, connected or non-connected. Since 
$Z_A ^\mathbb{C}(0) = 1$, 
the formal logarithm is well-defined for
$Z_A ^\mathbb{C}(t)$. The graphical expansion of 
$\log Z_A ^\mathbb{C}(t)$
then consists of connected ribbon graphs.

\begin{thm}
\label{thm:HAgraphexp}
The graphical expansion of the logarithm of the
generating function $Z_A ^\mathbb{C}(t)$ of (\ref{eq:HAgenerating})
is given by
\begin{equation}
\label{eq:HAgraphicalexp}
\begin{split}
\log Z_A ^\mathbb{C}(t) &=
\log \int_{\mathcal{H}_A} e^{-\frac{1}{2}
\langle x^2 \rangle}
e^{\sum_j \frac{t_j}{j}\langle x^j\rangle}
d\mu(x)\\
&=
\sum_{\substack{\Gamma \text{ connected}\\
\text{ribbon graph}}}
\frac{1}{|\Aut_R(\Gamma)|} 
A_{g(\Gamma),f(\Gamma)} ^{or}\;\prod_j t_j ^{v_j(\Gamma)}\;.
\end{split}
\end{equation}
\end{thm}

Recall the graph theoretic formulas (\ref{eq:vandefromvj}).
If we change $t_j$ to $\beta t_j$, then the graph expansion 
receives an extra factor of
$\beta^{v(\Gamma)}$ in the contribution 
from
 $\Gamma$. If we change $e^{-\frac{1}{2}\langle x^2\rangle}$
to $e^{-\frac{\alpha}{2}\langle x^2\rangle}$ with a positive real number
$\alpha$, then a change of variable $x\mapsto x/\sqrt{\alpha}$
produces a factor $\alpha^{-e(\Gamma)}$ to the
$\Gamma$-contribution. Therefore, 
\begin{equation}
\label{eq:alphabeta}
\log \int_{\mathcal{H}_A} e^{-\frac{\alpha}{2}
\langle x^2 \rangle}
e^{\beta\sum_j \frac{t_j}{j}\langle x^j\rangle}
d\mu_\alpha(x)
=
\sum_{\substack{\Gamma \text{ connected}\\
\text{ribbon graph}}}
\frac{\alpha^{-e(\Gamma)}\beta^{v(\Gamma)}}{|\Aut_R(\Gamma)|} 
A_{g(\Gamma),f(\Gamma)} ^{or}\;\prod_j t_j ^{v_j(\Gamma)}\;,
\end{equation}
where the normalized Lebesgue measure is adjusted for
$e^{-\frac{\alpha}{2}\langle x^2\rangle}$. 
An important example of 
Theorem~\ref{thm:HAgraphexp} is a 
hermitian matrix integral.  

\begin{ex}
\label{ex:hermitianmatrix}
Let us apply (\ref{eq:HAgraphicalexp}) to a complex 
matrix algebra
$A=M(n,\mathbb{C})$. 
The $*$-operation on this algebra is the matrix adjoint, and
\begin{equation}
\label{eq:GUEtrace}
\langle X \rangle = \frac{1}{n}\; \tr X
\end{equation}
is the normalized 
trace. 
The space of self-adjoint elements is the set of hermitian
matrices:
$$
\mathcal{H}_{M(n,\mathbb{C})} = \mathcal{H}_{n,\mathbb{C}}
\;.
$$
As an orthonormal basis, we use 
 $\{\sqrt{n}e_{ij}\}$, where
$$
e_{ij} = \big[\delta_{i\alpha} \delta_{j\beta}\big]_{\alpha,\beta}
$$
is the $n\times n$ elementary matrix which
has  $1$ at its $ij$ entry and $0$
everywhere else. 
Then we have
\begin{equation}
\label{eq:GUE}
\log \int_{\mathcal{H}_{n,\mathbb{C}}}
e^{-\frac{1}{2}\tr (X^2)}
e^{\sum_j \frac{t_j}{j} \tr (X^j)} d\mu(X)\\
=\sum_{\substack{\Gamma \text{ connected}\\
\text{ribbon graph}}}
\frac{1}{|\Aut_R(\Gamma)|} n^{f(\Gamma)}
\prod_j t_j ^{v_j(\Gamma)}\;.
\end{equation}
Indeed, the computation of $M(n,\mathbb{C})_{g,f} ^{or}$
is just evaluating the trace of the identity matrix
$I$.
Each tadpole contributes
$$
\sum_{i,j} e_{ij}e_{ij} ^* =  \sum_{i,j} e_{ii}=n\cdot I\;,
$$
while each bi-petal flower contributes
$$
\sum_{i,j,k,\ell} e_{ij} e_{k\ell} e_{ji} e_{\ell k}=
\sum_{i,j}e_{ii} e_{jj} = I\;.
$$
Therefore, we have 
$$
M(n,\mathbb{C})_{g,f} ^{or}= n^{2g} n^{2(f-1)}  \langle I\rangle
 = n^{-v+e+f}\;.
$$
Eqn.~(\ref{eq:GUE}) follows from (\ref{eq:GUEtrace}) and
(\ref{eq:alphabeta}).
Another useful form of hermitian matrix integral is
\begin{equation}
\label{eq:nGUE}
\log \int_{\mathcal{H}_{n,\mathbb{C}}}
e^{-\frac{n}{2}\tr (X^2)}
e^{n\sum_j \frac{t_j}{j} \tr (X^j)} d\mu(X)\\
=\sum_{\substack{\Gamma \text{ connected}\\
\text{ribbon graph}}}
\frac{1}{|\Aut_R(\Gamma)|} \;n^{\rchi(S_\Gamma)}
\prod_j t_j ^{v_j(\Gamma)}\;,
\end{equation}
which also follows from (\ref{eq:alphabeta}).
\end{ex}

Eqn.~(\ref{eq:GUE}) is due to \cite{BIZ} and has been used by many
authors in the study of hermitian matrix integrals
\cite{Harer-Zagier, Mulase1995, O, OP1, OP2, Penner, W1991a}.
In Section~\ref{sect:generating}, we give another 
example of the general formula, where we consider
$A=\mathbb{C}[G]$. 

\bigskip

\section{Integrals over a real von Neumann algebra}
\label{sect:realvna}

For a finite-dimensional real
von Neumann algebra $A$, the corresponding 
formulas become quite different. 
Since its trace is real valued, the hermitian inner product is real symmetric:
$$
\langle a,b\rangle = \langle ab^*\rangle
= \langle ba^* \rangle = \langle b,a\rangle\;.
$$
The integral we wish to evaluate is
\begin{equation}
\label{eq:realintegral}
\sum_{\substack{
(v_1,v_2,v_3,\cdots)\in\mathbb{N}^\infty\\
\mathbf{e}(v_1,v_2,v_3,\cdots)=n}}
\int_{\mathcal{H}_A} e^{-\frac{1}{4}\langle x^2\rangle}
\prod_{j\ge1} ^{\fin}
\frac{1}{v_j!\cdot (2j)^{v_j}}
 \langle x^j\rangle^{v_j}d\mu(x)
\end{equation}
for every $n\ge 0$
with respect to a different normalized Lebesgue measure
\begin{equation}
\label{eq:realnormal}
d\mu(x) = \frac{dx}{\int_{\mathcal{H}_A}
e^{-\frac{1}{4}\langle x^2\rangle} dx}\;.
\end{equation}
The generating function for (\ref{eq:realintegral})  is given by
\begin{equation}
\label{eq:ZRt}
\begin{split}
Z_A ^\mathbb{R}(t) 
&= \int_{\mathcal{H}_A} e^{-\frac{1}{4}\langle x^2\rangle}
e^{\sum_{j=1} ^\infty  \frac{t_j}{2j}\langle x^j\rangle}\;d\mu(x)\\
&=\sum_{n=0} ^\infty
\sum_{\substack{
(v_1,v_2,v_3,\cdots)\in\mathbb{N}^\infty\\
\mathbf{e}(v_1,v_2,v_3,\cdots)=n}}
\prod_{j\ge 1} ^\fin \frac{1}{v_j!\cdot (2j)^{v_j}}
\int_{\mathcal{H}_A} e^{-\frac{1}{4}\langle x^2\rangle}
\prod_{j\ge 1} ^\fin \langle x^j\rangle^{v_j}d\mu(x)\;.
\end{split}
\end{equation}

\begin{lem}
\label{lem:reallaplace}
Let $A$ be a real von Neumann algebra. Then
\begin{equation}
\label{eq:reallaplace}
\int_{\mathcal{H}_A} e^{-\frac{1}{4}\langle x^2\rangle}
\prod_{j=1} ^n \langle x^j\rangle^{v_j}d\mu(x)
=\left. \prod_{j=1} ^n 
\left\langle \left(
\frac{\partial}{\partial y}\right)^j\right\rangle^{v_j}
e^{\frac{1}{4}\langle(y+y^*)^2\rangle}\right|_{y=0}\;.
\end{equation}
\end{lem}

\begin{proof}
The adjoint of the element
$y\in A$ of Lemma~\ref{lem:exy} is given by
$$
y^* = \sum_{i-1} ^N y_i\; e_i ^*\;.
$$
We note here that $y_i\in\mathbb{R}$ and the orthonormal
basis $\{e_1,\dots,e_N\}$ for $A$ is a real basis. 
Unlike the complex case, we have
$\langle xy^*\rangle = \langle x^* y\rangle$, and hence
$$
\frac{\partial}{\partial y}\; e^{\langle xy^*\rangle}
=
\frac{\partial}{\partial y}\; e^{\langle x^*y\rangle}
=x^* e^{\langle xy^*\rangle}
$$
from Lemma~\ref{lem:exy}. Therefore, for $x\in\mathcal{H}_A$,
$$
\frac{\partial}{\partial y}\; e^{\langle x(y+y^*)/2\rangle}
=\frac{x+x^*}{2}e^{\langle x(y+y^*)/2\rangle}
= x\;e^{\langle x(y+y^*)/2\rangle}\;.
$$
The completion of the square is modified to 
\begin{equation*}
\begin{split}
\int_{\mathcal{H}_A} e^{-\frac{1}{4}\langle x^2\rangle}
e^{\langle x(y+y^*)/2\rangle} d\mu(x)
&=
\int_{\mathcal{H}_A} e^{-\frac{1}{4}\langle (x-(y+y^*))^2\rangle}
e^{\frac{1}{4}\langle (y+y^*)^2\rangle} d\mu(x)\\
&=e^{\frac{1}{4}\langle (y+y^*)^2\rangle}\;.
\end{split}
\end{equation*}
The rest of the proof is the same as the complex case.
\end{proof}

\begin{figure}[htb]
\centerline{\epsfig{file=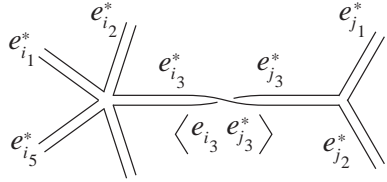, width=2in}}
\caption{A half-edge labeled by $i_3$ of the left vertex
is connected to a half-edge labeled by $j_3$ of the
vertex at the right by a twisted edge. A propagator
$\langle e_{i_3} e_{j_3} ^*\rangle$ is assigned to this edge.}
\label{fig:twistedge}
\end{figure}

To compute the RHS of (\ref{eq:reallaplace}), 
we first note
\begin{equation*}
\begin{split}
\frac{\partial}{\partial y} \;e^{\frac{1}{4}\langle
(y+y^*)^2\rangle} &= \frac{1}{2}
\sum_i \frac{\partial}{\partial y_i}\; e_i ^*\;
e^{\frac{1}{4}\langle (\sum_j y_j\; e_j + \sum_j 
{y_j}\; e_j ^*)^2
\rangle}\\
&= \frac{1}{2}
\sum_i \langle (y+y^*) (e_i+e_i ^*) \rangle \; e_i ^*\;
e^{\frac{1}{4}\langle(y+y^*)^2\rangle}\\
&= \frac{1}{2}
\sum_i \langle (y+y^*) ,(e_i+e_i ^*) \rangle \; e_i ^*\;
e^{\frac{1}{4}\langle(y+y^*)^2\rangle}\\
&=(y+y^*) \; e^{\frac{1}{4}\langle(y+y^*)^2\rangle}\;.
\end{split}
\end{equation*}
In particular, 
\begin{equation*}
\frac{\partial}{\partial y_i}
e^{\frac{1}{4}\langle(y+y^*)^2\rangle}
= \frac{1}{2}\; \langle (y+y^*) ,(e_i+e_i ^*) \rangle \;
e^{\frac{1}{4}\langle(y+y^*)^2\rangle}\;,
\end{equation*}
and hence
\begin{equation}
\label{eq:dyidyjreal}
\begin{split}
\left.\frac{\partial}{\partial y_i}\frac{\partial}{\partial y_j}
e^{\frac{1}{4}\langle(y+y^*)^2\rangle}\right|_{y=0}
&=  \frac{1}{2}
\left.\frac{\partial}{\partial y_i}
\langle (y+y^*) ,(e_j+e_j ^*) \rangle 
\right|_{y=0}\\
&= \frac{1}{2}\;\langle (e_i+e_i ^*) ,
e_j+e_j ^*)\rangle\\
&=\langle e_i e_j ^*\rangle +\langle e_i e_j \rangle\\
&= \delta_{ij}+\langle e_i e_j \rangle\;.
\end{split}
\end{equation}
This formula has an extra term $\delta_{ij}$ compared to
(\ref{eq:dyidyjcomplex}). 
To compute the graphical expansion of 
\begin{equation}
\label{eq:realcompute}
\sum_{\substack{
(v_1,v_2,v_3,\cdots)\in\mathbb{N}^\infty\\
\mathbf{e}(v_1,v_2,v_3,\cdots)=n}}
\prod_{j\ge 1} ^\fin \frac{1}{v_j!\cdot (2j)^{v_j}}
\left. 
\left\langle \left(
\frac{\partial}{\partial y}\right)^j\right\rangle^{v_j}
e^{\frac{1}{4}\langle(y+y^*)^2\rangle}\right|_{y=0}\;,
\end{equation}
we proceed as before and assign a cyclically ordered $j$-valent
vertex
to each factor $\langle (\partial/\partial y)^j\rangle$
of the differentiation,
assign  the vertex contribution
$\langle e_{i_1} ^*\cdots  e_{i_j} ^*\rangle$ to it, 
   and place the vertex
consistently on an oriented plane with the clockwise orientation. 
When a pair of differentiation is applied, because of 
(\ref{eq:dyidyjreal}), there are now two choices:
straight
connection as in the ribbon graph case
Figure~\ref{fig:edgeprop} with a propagator
$\langle e_i e_j\rangle$ assigned to the edge,
 or connection with a twisted edge
carrying a propagator $\langle e_i e_j ^*\rangle = \delta_{ij}$
as in Figure~\ref{fig:twistedge}.

If the straight connection is used at every edge, then the resulting
graph is a ribbon graph as before. Otherwise, we obtain a
M\"obius graph $\Gamma$ with some twisted edges. 
It has to be noted that the existence of  twisted edges does not
necessarily mean that the M\"obius graph is non-orientable. 
For a M\"obius graph $\Gamma$ thus obtained, let us define
the graph contribution $A_\Gamma ^\mathbb{R}$ as
the sum with respect to all indices of the product of 
all vertex contributions $\langle e_{i_1} ^*\cdots  e_{i_j} ^*\rangle$
and the product of all propagators, where 
$\langle e_ie_j\rangle$ is chosen for a straight edge and
$\langle e_ie_j ^*\rangle$ is chosen for a twisted edge.

The reality condition of the trace provides another invariance:
\begin{equation}
\label{eq:vertexreality}
\left\langle e_{i_1} ^* e_{i_2} ^*\cdots e_{i_{j-1}} ^*
e_{i_j} ^*\right\rangle
=
\left\langle e_{i_j}  e_{i_{j-1}} \cdots
e_{i_2} e_{i_1}\right\rangle\;.
\end{equation}
This equality brings an equivalence relation into the
set of M\"obius  graphs.
To identify it, let us observe the following:

\begin{lem}
\label{lem:vertexflip}
Let us denote by $\langle e_i e_j ^{\pm *}\rangle$ either
$\langle e_ie_j\rangle$ or $\langle e_ie_j ^*\rangle$,
and use $\langle e_i e_j ^{\mp *}\rangle$ to indicate the other 
propagator. Then
\begin{equation}
\label{eq:vertexflip}
\sum_{i_1,\dots,i_j} \left\langle
e_{i_1} ^*\cdots e_{i_j} ^*\right\rangle
\langle e_{i_1} e_{h_1} ^{\pm *}\rangle\cdots
\langle e_{i_j} e_{h_j} ^{\pm *}\rangle
= \sum_{i_1,\dots,i_j} \left\langle
e_{i_j} ^*\cdots e_{i_1} ^*\right\rangle
\langle e_{i_1} e_{h_1} ^{\mp *}\rangle\cdots
\langle e_{i_j} e_{h_j} ^{\mp *}\rangle\;.
\end{equation}
\end{lem}

\begin{proof}
Using the contraction formula (\ref{eq:ab}), the LHS is equal to
$\left\langle
e_{h_1} ^{\pm *}\cdots e_{h_j} ^{\pm *}\right\rangle$.
Similarly, the RHS is equal to 
$\left\langle
e_{h_j} ^{\mp *}\cdots e_{h_1} ^{\mp *}\right\rangle$.
Because of the reality condition (\ref{eq:vertexreality}), 
these are actually the same.
\end{proof}

Notice that the equation  (\ref{eq:vertexflip}) is
exactly the vertex flip operation of 
Figure~\ref{fig:vertexflip}. 
This allows us to define the graph contribution $A_\Gamma ^\mathbb{R}$
slightly differently: it is the sum 
with respect to all indices of the 
product of all vertex
contributions $\langle e_{i_1} ^*\cdots  e_{i_j} ^*\rangle$
each of which has a cyclic order that is determined
according to the cyclic order of the vertex, and the product of
all propagators of edges determined by their twist. 
The extra redundancy of the vertex flip is compensated with 
the factor $(2j)^{v_j}$ in front of (\ref{eq:realcompute}),
which is the order of the product of dihedral groups acting
on the vertices  through rotations and flips. In a parallel way
with Proposition~\ref{prop:ribbonexp}, we have thus established:

\begin{prop}
\label{prop:moebexp}
Let $MG(e)$ be the set of all M\"obius
graphs consisting of $e$ edges.
For each $j$-valent vertex of  $\Gamma$, let
us assign 
$\big\langle e_{i_1} ^*\cdots e_{i_j} ^*\big\rangle$,
where the cyclic order of the product is determined by
the cyclic order of the vertex.
For every edge we assign 
a propagator $\langle e_{i_k}e_{h_\ell}\rangle$
if the edge is straight and 
 $\langle e_{i_k}e_{h_\ell} ^*\rangle$ if it is twisted.
The incidence relation should be consistent with the
labeling of half-edges, namely, the edge labeled with
$i_k h_\ell$
connects the half-edge labeled by $i_k$ of a vertex
to the half-edge labeled by $h_\ell$ from another
vertex.
Let $A_\Gamma ^{\mathbb{R}}$ 
denote the sum 
with respect to all indices of the product
of all contributions from vertices and edges. Then 
\begin{equation}
\label{eq:graphcontributionreal}
\begin{split}
&\quad \sum_{\substack{
(v_1,v_2,v_3,\cdots)\in\mathbb{N}^\infty\\
\mathbf{e}(v_1,v_2,v_3,\cdots)=2e}}
\int_{\mathcal{H}_A} e^{-\frac{1}{4}\langle x^2\rangle}
\prod_{j\ge 1} ^\fin 
\frac{\langle x^j\rangle^{v_j}}{v_j !\cdot (2 j)^{v_j}}d\mu(x)\\
&= 
\sum_{\substack{
(v_1,v_2,v_3,\cdots)\in\mathbb{N}^\infty\\
\mathbf{e}(v_1,v_2,v_3,\cdots)=2e}}
\prod_{j\ge 1} ^\fin \frac{1}{v_j !\cdot (2j)^{v_j}}
\left.\left\langle \left(
\frac{\partial}{\partial y}\right)^j\right\rangle^{v_j}
e^{\frac{1}{4}\langle(y+y^*)^2\rangle}\right|_{y=0}\\
&= \sum_{\Gamma\in MG(e)}
\frac{1}{|\Aut_M (\Gamma)|} \;A_\Gamma ^{\mathbb{R}}\;,
\end{split}
\end{equation}
where $\Aut_M (\Gamma)$
is the automorphism group of $\Gamma$ as a M\"obius graph. 
\end{prop}

If $\Gamma$ is orientable, then a series of vertex flip 
operations makes $\Gamma$ a ribbon graph, and for
such a graph, $A_\Gamma ^\mathbb{R} = A_\Gamma ^{or}$.
Although the von Neumann algebra $A$ is real, we can use the
same definition of $A_\Gamma ^{or}$ 
as in Proposition~\ref{prop:ribbonexp} for a real $A$. 
Its invariance with respect to the topological type of the
orientable surface $S_\Gamma$ is the same as before. 
Even a M\"obius graph
$\Gamma$ is non-orientable, we still have the following:

\begin{lem}
\label{lem:moebtopology}
Let $\Gamma$ be a connected 
M\"obius  graph with two or
more vertices, and $E$ an edge of $\Gamma$
incident to two distinct vertices. Then the contribution
of the graph $A_\Gamma ^{\mathbb{R}}$ 
is invariant under the edge-contraction:
\begin{equation}
\label{eq:moebtopology}
A_\Gamma ^{\mathbb{R}}= 
A_{\Gamma/ E} ^{\mathbb{R}}\;.
\end{equation}
\end{lem}

\begin{proof}
Let $V_1$ and $V_2$ be the two vertices incident
to the edge $E$. If $E$ is not twisted, then the same
argument of Lemma~\ref{lem:topology} applies. 
If the edge is twisted, then first apply a vertex flip 
operation to $V_2$ and untwist $E$. Then the
situation is the same as before, and we can contract
the edge,  joining $V_1$ and $V_2$ together.
We give the cyclic oder of $V_1$
as the new cyclic order to the newly created vertex.
\end{proof}

\begin{figure}[htb]
\centerline{\epsfig{file=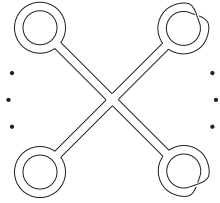, width=1.2in}}
\caption{A standard graph for a non-orientable surface
of topological type $(k,f)$. It has $f-1$ tadpoles on the left,
and $k$ twisted tadpoles on the right.}
\label{fig:standardnor}
\end{figure}

It is known \cite{Hatcher} that the set of all
connected M\"obius graphs with $f$ faces
drawn on a closed non-orientable
surface of cross-cap genus $k$ is \emph{connected}
with respect to the edge contraction and edge expansion 
moves. (These moves are called
\emph{diagonal flips} in \cite{Hatcher}.)
Therefore, we can compute the invariant
$A_\Gamma ^{\mathbb{R}}$ for a non-orientable
M\"obius graph again by choosing our favorite graph. 
If we use a graph of Figure~\ref{fig:standardnor}, then
for every non-orientable M\"obius graph of topological
type $(k,f)$, the graph contribution is equal to
\begin{equation}
\label{eq:Akf}
A_{k,f} ^{nor}
= \sum_{\substack{i_1,\dots,i_k\\
h_1,\dots, h_{f-1}}}
\left\langle
(e_{i_1} ^*)^2\cdots (e_{i_k} ^*)^2\cdot
e_{h_1} e_{h_1} ^* \cdots e_{h_{f-1}} e_{h_{f-1}} ^* 
\right\rangle\;.
\end{equation}

Now we have 

\begin{thm}
\label{thm:HArealexp}
The graphical expansion of the logarithm of the
generating function $Z_A ^\mathbb{R}(t)$ of 
(\ref{eq:ZRt}) associated with a real von Neumann 
algebra $A$
is given by
\begin{equation}
\label{eq:HArealexp}
\begin{split}
\log Z_A ^\mathbb{R}(t) &=
\log \int_{\mathcal{H}_A} e^{-\frac{1}{4}
\langle x^2 \rangle}
e^{\sum_j \frac{t_j}{2j}\langle x^j\rangle}
d\mu(x)\\
&=
\sum_{\substack{\Gamma \text{ connected orientable}\\
\text{M\"obius graph}}}
\frac{1}{|\Aut_M(\Gamma)|} 
A_{g(\Gamma),f(\Gamma)} ^{or}
\;\prod_j t_j ^{v_j(\Gamma)}\\
&+
\sum_{\substack{\Gamma \text{ connected non-}\\
\text{orientable M\"obius graph}}}
\frac{1}{|\Aut_M(\Gamma)|} 
A_{k(\Gamma),f(\Gamma)} ^{nor}
\; \prod_j t_j ^{v_j(\Gamma)}\;.
\end{split}
\end{equation}
\end{thm}

For future convenience, 
we also record 
\begin{equation}
\label{eq:realalphabeta}
\begin{split}
&
\log \int_{\mathcal{H}_A} e^{-\frac{\alpha}{4}
\langle x^2 \rangle}
e^{\beta\sum_j \frac{t_j}{2j}\langle x^j\rangle}
d\mu(x)\\
&=
\sum_{\substack{\Gamma \text{ connected orientable}\\
\text{M\"obius graph}}}
\frac{\alpha^{-e(\Gamma)}\beta^{v(\Gamma)}}{|\Aut_M(\Gamma)|} 
A_{g(\Gamma),f(\Gamma)} ^{or}
\;\prod_j t_j ^{v_j(\Gamma)}\\
&+
\sum_{\substack{\Gamma \text{ connected non-}\\
\text{orientable M\"obius graph}}}
\frac{\alpha^{-e(\Gamma)}\beta^{v(\Gamma)}}{|\Aut_M(\Gamma)|} 
A_{k(\Gamma),f(\Gamma)} ^{nor}
\; \prod_j t_j ^{v_j(\Gamma)}
\end{split}
\end{equation}
for a positive $\alpha$ and any $\beta$.

As an example of these formulas, let us consider
the case when $A$ is a simple algebra. 
This time, it is isomorphic to 
either $M(n,\mathbb{R})$ or $M(n,\mathbb{H})$. 

\begin{ex}
\label{ex:realsym}
Let $A=M(n,\mathbb{R})$. Then as in 
Example~\ref{ex:hermitianmatrix}, we can use
$\{\sqrt{n}e_{ij}\}$ as our orthonormal basis.
 The space $\mathcal{H}_A$
of self-adjoint elements is the set of all real
symmetric matrices
$\mathcal{H}_{n,\mathbb{R}}$.
 Since elementary matrices are defined
over the reals, we immediately see
$$
A_{g,f} ^{or} = n^{-v+e+f}
$$
as before. To calculate $A_{k,f} ^{nor}$, we note that
$$
(e_{i_1 j_1})^2\cdots (e_{i_k j_k})^2 = 0
$$
unless all $2k$ indices are the same, and if they are the same, 
then the result is $e_{ii}$. Thus the sum of all these products
is the identity matrix $I$. The contribution from 
the tadpoles of Figure~\ref{fig:standardnor} is the same
as in Example~\ref{ex:hermitianmatrix}, so we have
\begin{equation}
\label{eq:GOE}
\log \int_{\mathcal{H}_{n,\mathbb{R}}}
e^{-\frac{1}{4}\tr (X^2)}
e^{\sum_j \frac{t_j}{2j} \tr (X^j)} d\mu(X)\\
=\sum_{\substack{\Gamma \text{ connected}\\
\text{M\"obius graph}}}
\frac{1}{|\Aut_M(\Gamma)|} \;n^{f(\Gamma)}\;
\prod_j t_j ^{v_j(\Gamma)}\;,
\end{equation}
or more conveniently, 
\begin{equation}
\label{eq:nGOE}
\log \int_{\mathcal{H}_{n,\mathbb{R}}}
e^{-\frac{n}{4}\tr (X^2)}
e^{n\sum_j \frac{t_j}{2j} \tr (X^j)} d\mu(X)\\
=\sum_{\substack{\Gamma \text{ connected}\\
\text{M\"obius graph}}}
\frac{1}{|\Aut_M(\Gamma)|} \;n^{\rchi(S_\Gamma)}\;
\prod_j t_j ^{v_j(\Gamma)}\;.
\end{equation}
These formulas are well known  (see for example,
\cite{BIPZ, GHJ}). 
\end{ex}

\begin{ex}
\label{ex:quatselfadjoint}
This time let us choose $A=M(n,\mathbb{H})$. 
As a real basis for quaternions, we use
\begin{equation}
\label{eq:quatbasis}
e^0 = \begin{pmatrix}
1\\
&1
\end{pmatrix}\;,\quad
e^1 = \begin{pmatrix}
&1\\
-1
\end{pmatrix}\;,\quad
e^2 = \begin{pmatrix}
i\\
&-i
\end{pmatrix}\;,\quad
e^3 = \begin{pmatrix}
&-i\\
-i
\end{pmatrix}\;.
\end{equation}
The adjoint operation $(e^\mu)^*$ on these $2\times 2$ 
matrices is the same
as the
 conjugate transposition $(e^\mu)^\dagger$. 
A real basis for $M(n,\mathbb{H})$
is given by
$$
e_{ij} ^\nu = e_{ij}\tensor e^\nu, \qquad i,j = 1, \dots, n;
\quad\nu = 0, \dots,3\;.
$$
The normalized trace is defined by 
$$
\langle e_{ij} ^\nu \rangle = \frac{1}{2n}\; \tr_{n\times n} (e_{ij})
\cdot \tr_{2\times 2}( e^\nu)
$$
for the basis elements and $\mathbb{R}$-linearly extended to 
all matrices. We notice that $\langle\; \rangle$
 is real valued because
the $2\times 2$ trace has value $0$ for imaginary quaternionic units.
With respect to the normalized trace, $\{\sqrt{n}e_{ij} ^\nu\}$ is an orthonormal basis. The $*$-operation with respect the basis
is given by
$$
(e_{ij} ^\nu)^* = e_{ji}\tensor (e^\nu)^\dagger\;.
$$
The space of self-adjoint elements 
$\mathcal{H}_A = \mathcal{H}_{n,\mathbb{H}}$ consists
of self-adjoint quatermionic matrices of size $n\times n$, and is
spanned by $e_{ij}\tensor e^\nu + e_{ji}\tensor (e^\nu)^\dagger$. 
Note that the diagonal entries of a self-adjoint matrix are spanned by
$e^\nu + (e^\nu)^\dagger$, and hence are real. Thus we have a real
linear map
$$
\tr_{n\times n} :\mathcal{H}_{n,\mathbb{H}}\longrightarrow
\mathbb{R}\;.
$$
Since $e^\nu(e^\nu)^\dagger = e^0$, we have
$$
\sum_{i,j,\nu} e_{ij} ^\nu (e_{ij} ^\nu)^* =4nI\tensor e^0\;.
$$
Similarly, $e^0 = (e^0)^2= -(e^1)^2=-(e^2)^2=-(e^3)^2$, 
hence
$$
\sum_{i,j,\nu}( e_{ij} ^\nu)^2 = -2 I\tensor e^0\;.
$$
To compute $e_{ij} ^\mu e_{k\ell} ^\nu (e_{ij} ^\mu)^*
( e_{k\ell} ^\nu)^*$, we use
\begin{equation}
\label{eq:munumunu}
e^\mu e^\nu (e^\mu)^\dagger (e^\nu)^\dagger
=\begin{cases} 
-e^0\qquad \;\mu,\nu>0, \; \mu\ne \nu\\
e^0\qquad\quad \text{otherwise}\;.
\end{cases}
\end{equation}
Therefore, of the $16$ combinations, $6$ cases are equal to 
$-e^0$ and $10$ are equal to $e^0$. Thus
$$
\sum_{\mu,\nu} e^\mu e^\nu (e^\mu)^\dagger (e^\nu)^\dagger
=4e^0\;,
$$
and altogether, we have
$$
\sum_{i,j,k,\ell,\mu,\nu} e_{ij} ^\mu e_{k\ell} ^\nu (e_{ij} ^\mu)^*
( e_{k\ell} ^\nu)^*
=4I\tensor e^0\;.
$$
From all the above, we calculate
\begin{equation*}
\begin{split}
M(n,\mathbb{H})_{g,f} ^{or}
&= (4n^2)^g (4n^2)^{f-1}\langle I\tensor e^0\rangle
= (2n)^{-v+e+f}\\
M(n,\mathbb{H})_{k,f} ^{nor}
&= (-2n)^k (4n^2)^{f-1}\langle I\tensor e^0\rangle
= (-1)^k (2n)^{-v+e+f}\;.
\end{split}
\end{equation*}
Note that $(-1)^k = (-1)^{\rchi(S_\Gamma)}$. 
Combining these computations
with (\ref{eq:realalphabeta}) and using the
$n\times n$ trace
of $M(n,\mathbb{H})$, we finally obtain
\begin{equation}
\label{eq:GSE}
\log \int_{\mathcal{H}_{n,\mathbb{H}}}
e^{-\frac{1}{2}\tr (X^2)}
e^{\sum_j \frac{t_j}{j} \tr (X^j)} d\mu(X)\\
=\sum_{\substack{\Gamma \text{ connected}\\
\text{M\"obius graph}}}
\frac{(-1)^{\rchi(S_\Gamma)}}{|\Aut_M(\Gamma)|} \;(2n)^{f(\Gamma)}\;
\prod_j t_j ^{v_j(\Gamma)}\;,
\end{equation}
or equivalently, 
\begin{equation}
\label{eq:nGSE}
\log \int_{\mathcal{H}_{n,\mathbb{H}}}
e^{-n\tr (X^2)}
e^{2n\sum_j \frac{t_j}{j} \tr (X^j)} d\mu(X)\\
=\sum_{\substack{\Gamma \text{ connected}\\
\text{M\"obius graph}}}
\frac{1}{|\Aut_M(\Gamma)|} \;(-2n)^{\rchi(S_\Gamma)}\;
\prod_j t_j ^{v_j(\Gamma)}\;.
\end{equation}
\end{ex}

These results are in agreement with recently established
formulas found in \cite{MW}. 

\bigskip

\section{Generating functions for the number of
representations of surface groups}
\label{sect:generating}

Let us now turn our attention to the case of the complex
group algebra $A=\mathbb{C}[G]$
of a finite group $G$. The $*$-operation is defined by
\begin{equation}
\label{eq:realstar}
*:\mathbb{C}[G]\owns x=\sum_{w \in G} x(w)  \cdot
w  \longmapsto x^* = \sum_{w \in G} 
\overline{x(w) } \cdot
w ^{-1}\in \mathbb{C}[G].
\end{equation}
As the trace, we use
\begin{equation}
\label{eq:Ctrace}
\langle\; \rangle = \frac{1}{|G|} \rchi_{\reg},
\end{equation}
where $\rchi_\reg$ is the character of the regular representation
of $G$ on $\mathbb{C}[G]$, linearly extended  to the
whole group algebra. 
The self-adjoint condition $x^*=x$ means
$x(w ^{-1})=\overline{x(w) }$, and 
 we have 
$\mathcal{H}_{\mathbb{C}[G]} = \mathbb{R}^{|G|}$
as a real vector space. A natural orthonormal basis
for $\mathbb{C}[G]$ is the group $G$ itself, since we have
\begin{equation*}
\langle uv^*\rangle = \frac{1}{|G|}\rchi_\reg(uv^{-1})
=
\begin{cases}
1\qquad u=v\\
0\qquad\text{otherwise}\;.
\end{cases}
\end{equation*}
It is because the normalized trace on $G$ takes value $1$ only when
the group element is the identity and $0$ otherwise. 
Recall that 
$$
\pi_1(S) = \big\langle a_1,b_1,\dots,a_g,b_g\;\big|\;
a_1b_1a_1 ^{-1}b_1 ^{-1}\cdots
a_gb_ga_g ^{-1}b_g ^{-1} = 1\big\rangle\;,
$$
where $S$ is an orientable surface of genus $g$.
With respect the orthonormal basis, we immediately see
\begin{equation}
\label{eq:CGgf}
\begin{split}
\mathbb{C}[G]_{g,f} ^{or}&=
\sum_{u_i,v_i,w_j\in G}
\left\langle u_1v_1u_1 ^{-1}v_1 ^{-1}\cdots
u_gv_gu_g ^{-1}v_g ^{-1}\cdot
w_1w_1 ^{-1}\cdots w_{f-1} w_{f-1} ^{-1}\right\rangle\\
&=|G|^{f-1}\;|\Hom(\pi_1(S),G)|\;.
\end{split}
\end{equation}
Summarizing these facts and changing the 
constant factors using
(\ref{eq:alphabeta}), we obtain a
generating function of the number of homomorphisms
from the fundamental group of an orientable surface
into the finite group.

\begin{thm}
\label{thm:CGexpansion}
Let $G$ be a finite group. 
The following integral over the self-adjoint elements of 
the complex group algebra $\mathbb{C}[G]$ gives
the generating function for the cardinality of the
representation variety of an orientable surface group 
in $G$:
\begin{multline}
\label{eq:CGasymptotic}
\log \int_{\mathcal{H}_{\mathbb{C}[G]}}
\exp\left(-\frac{1}{2}\;\rchi_\reg(  x^2 )\right)
\exp\left(\sum_{j} \frac{t_j}{j}\rchi_\reg(x^j )
\right) d\mu(x)\\
= \sum_{\substack{\Gamma  \text{ connected}\\
\text{ribbon graph}}}
\frac{1}{|\Aut_R \Gamma|}|G|^{\rchi(S_\Gamma)-1}
|\Hom(\pi_1(S_\Gamma),G)|
\prod_{j} t_j ^{v_j(\Gamma)}\;.
\end{multline}
\end{thm}

Note that we have a von Neumann algebra isomorphism
\begin{equation}
\label{eq:complexdecomp}
\mathbb{C}[G] \isom \bigoplus_{\lambda\in\hat{G}}
\End(V_\lambda)\;,
\end{equation}
which decomposes the character of the regular representation
into the sum of irreducible characters: 
$$
\rchi_{\reg} = \sum_{\lambda\in\hat{G}} 
(\dim \lambda) \;\rchi_\lambda 
= \sum_{\lambda\in\hat{G}} (\dim \lambda) \;
\tr_{\lambda}\;,
$$
where $ \dim \lambda$ is the dimension of $\lam\in\Ghat$ and 
$\rchi_\lambda$ is its character.
Therefore, using (\ref{eq:nGUE}) for each 
irreducible factor, we calculate
\begin{equation}
\label{eq:oricomputation}
\begin{split}
&\log \int_{\mathcal{H}_{\mathbb{C}[G]}}
\exp\left(-\frac{1}{2}\;\rchi_\reg(  x^2 )\right)
\exp\left(\sum_{j} \frac{t_j}{j}\rchi_\reg (x^j )
\right) d\mu(x)\\
=&\log \int_{\mathcal{H}_{\mathbb{C}[G]}}
\prod_{\lambda\in\hat{G}}
\exp\left(-\frac{\dim \lam}{2}\;\tr_{\lambda}(  x^2 )\right)
\exp\left(\dim \lam \sum_{j} 
\frac{t_j}{j}\tr_{\lambda}(x^j )
\right) d\mu_\lambda(x)\\
=&\sum_{\lambda\in\hat{G}}\log 
\int_{\mathcal{H}_{\dim \lam, \mathbb{C}}}
\exp\left(-\frac{\dim \lam}{2}\;\tr_{\lambda}(  x^2 )\right)
\exp\left(\dim \lam \sum_{j} 
\frac{t_j}{j}\tr_{\lambda}(x^j )
\right) d\mu_\lambda(x)\\
=&\sum_{\substack{\Gamma  \text{ connected}\\
\text{ribbon graph }}}
\frac{1}{|\Aut_R \Gamma|}
\sum_{\lambda\in\hat{G}}
 (\dim \lam)^{\rchi (S_\Gamma)}
\prod_{j} t_j ^{v_j(\Gamma)}\;,
\end{split}
\end{equation}
where $d\mu_\lambda$ is the normalized Lebesgue
measure on the space of
$\dim \lam\times \dim \lam$ hermitian matrices.
Comparing the two expressions (\ref{eq:CGasymptotic})
and (\ref{eq:oricomputation}),
we recover Mednykh's formula (\ref{eq:Med}):
$$
\sum_{\lam\in\Ghat} (\dim \lam)^{\rchi(S)}
= |G|^{\rchi(S)-1}|\Hom(\pi_1(S), G)|\;.
$$

\begin{rem}
Another proof of Mednykh's formula is found
 in \cite{FQ}, which uses Chern-Simons
gauge theory with a finite gauge group. 
Burnside asked a related question
on p.~319 (\S~238, Ex.~7) of 
his textbook \cite{Burnside}. The formula
for genus $1$ case is found
in Frobenius \cite{F} of 1896. We refer to  \cite{Stanley}
for the relation of these formulas to combinatorics. 
An excellent historical account on this and 
Frobenius-Schur formula (\ref{eq:FScount})
is found in \cite{Jones}.
\end{rem}

Now consider the real group algebra $\mathbb{R}[G]$. 
For a non-orientable surface of cross-cap genus $k$, we 
know
$$
\pi_1(S) = \big\langle
a_1,\dots,a_k\;\big|\;a_1 ^2\cdots a_k ^2 = 1\big\rangle\;.
$$
Therefore, 
\begin{equation}
\label{eq:RGkf}
\begin{split}
\mathbb{R}[G]_{k,f} ^{nor}
&= \sum_{u_i,w_j\in G}
\left\langle u_1 ^2\cdots u_k ^2\cdot w_1 ^{-1} w_1\cdots
w_{f-1} ^{-1} w_{f-1} \right\rangle\\
&= |G|^{f-1}\;|\Hom(\pi_1(S), G)|\;.
\end{split}
\end{equation}
Our general formula (\ref{eq:realalphabeta}) yields

\begin{thm}
Let $G$ be a finite group. The following integral
over the space of self-adjoint elements of the
real group algebra $\mathbb{R}[G]$ gives the
generating function for the number of homomorphisms
from the fundamental group of a closed surface into $G$,
$|\Hom(\pi_1(S), G)|$, for all $S$, including orientble
and non-orientable surfaces.
\begin{multline}
\label{eq:realgr}
\log  \int_{\mathcal{H}_{\mathbb{R}[G]}}
e^{-\frac{1}{4}\rchi_{\reg}( x^2)} 
e^{\frac{1}{2}\sum_j \frac{t_j}{j}
\rchi_{\reg}( x^j)}d\mu(x)\\
=
\sum_{\substack{\Gamma  \text{ connected}\\
\text{M\"obius graph}}}
\frac{1}{|\Aut_M \Gamma|} |G|^{\rchi(S_\Gamma)-1} 
|\Hom(\pi_1(S_\Gamma),G)|
\prod_{j} t_j ^{v_j (\Gamma)}\;.
\end{multline}
\end{thm}

Recall that the real group algebra $\mathbb{R}[G]$
decomposes into simple factors
according to the three types of irreducible representations
(\ref{eq:FSindicator}). 
Notice that $\hat{G}_1$ consists of complex irreducible
representations of $G$ that are defined over $\mathbb{R}$.
A representation in $\hat{G}_2$ is not defined over
$\mathbb{R}$, and its character is not real-valued. Thus
the complex conjugation acts on the set $\hat{G}_2$
without fixed points. Let $\hat{G}_{2+}$ denote a half
of $\hat{G}_2$ such that
\begin{equation}
\label{eq:G2+}
\hat{G}_{2+}\cup \overline{\hat{G}_{2+}} = \hat{G}_2\;.
\end{equation}
A complex irreducible 
representation of $G$ that belongs to $\hat{G}_4$
admits a skew-symmetric 
bilinear form. In particular, its dimension (over $\mathbb{C}$)
is even. Now we have a von Neumann algebra isomorphism
\begin{equation}
\label{eq:realdecomp}
\mathbb{R}[G] \isom \bigoplus_{\lambda\in \hat{G}_1}
\End_{\mathbb{R}}(\lam^\mathbb{R})
\dsum \bigoplus_{\lambda\in \hat{G}_{2+}}
\End_{\mathbb{C}}(\lambda)
\dsum \bigoplus_{\lambda\in \hat{G}_4}
\End_{\mathbb{H}}(\lambda ^\mathbb{H})\;,
\end{equation}
where $\lambda ^\mathbb{R}$ is a real irreducible
representation of $G$ that satisfies 
$\lambda =\lambda ^\mathbb{R}\tensor_\mathbb{R}
\mathbb{C}$. The representation
space $\lambda ^\mathbb{H}$ is a 
$(\dim \lambda)/2$-dimensional vector
space defined over $\mathbb{H}$ for
$\lambda\in\hat{G}_4$ such that its image 
under the natural injection
\begin{equation}
\label{eq:HtoC}
\End_\mathbb{H} (\lambda ^\mathbb{H})
\longrightarrow \End_\mathbb{C} (\lambda)
\end{equation}
coincides with the image of 
$$
\rho_\lambda : \mathbb{R}[G]\longrightarrow
\End_\mathbb{C} (\lambda)\;,
$$
where $\rho_\lam$ is the representation 
 of $\mathbb{R}[G]$ corresponding to
$\lam\in\Ghat$. The injective algebra homomorphism
(\ref{eq:HtoC}) is defined by the $2\times 2$ matrix
representation of the quaternions (\ref{eq:quatbasis}).
The algebra isomorphism (\ref{eq:realdecomp})
gives a formula for the character
of the regular representation on $\mathbb{R}[G]$:
\begin{equation}
\label{eq:realreg}
\rchi_{\reg} = \sum_{\lambda\in\hat{G}_1} (\dim
\lambda ) \;\rchi_\lambda +
 \sum_{\lambda\in\hat{G}_{2+}} (\dim
\lambda)\; (\rchi_\lambda +  \overline{\rchi_\lambda})+
\sum_{\lambda\in\hat{G}_4}2 (\dim 
\lambda) \cdot \trace_{\lambda ^\mathbb{H}}\; ,
\end{equation}
where in the last term the character 
is given as the trace of quaternionic
$(\dim\lambda)/2\times 
(\dim\lambda)/2$ matrices.
Notice that if  $\lam\in \Ghat_2$, then for every
$x = x^*\in \mathcal{H}_{\mathbb{R}[G]}$, we have
$$
\rchi_\lam (x) = \overline{\rchi_\lam}(x)
=\tr_{\dim\lam}(\rho_\lam(x))
$$
since $\rho_\lam(x)$ is a hermitian matrix of size $\dim\lam\times
\dim\lam$. 

The integration (\ref{eq:realgr}) can be carried out
using (\ref{eq:nGUE}), (\ref{eq:nGOE}) and (\ref{eq:nGSE})
with the decomposition (\ref{eq:realdecomp}) and 
(\ref{eq:realreg}). The result is
\begin{equation}
\label{eq:realgadecomp}
\begin{split}
\log  &\int_{\mathcal{H}_{\mathbb{R}[G]}}
e^{-\frac{1}{4}\rchi_{\reg}( x^2)} 
e^{\frac{1}{2}\sum_j \frac{t_j}{j}
\rchi_{\reg}( x^j)}d\mu(x)\\
&=\sum_{\lam\in\Ghat_1}
\log \int_{\mathcal{H}_{\dim\lam, \mathbb{R}}}
e^{-\frac{\dim\lam}{4}\tr(x^2)}
e^{\dim\lam \sum_j\frac{t_j}{2j}\tr(x^j)}d\mu_\lam(x)\\
&\quad+\sum_{\lam\in\Ghat_2}
\log \int_{\mathcal{H}_{\dim\lam, \mathbb{C}}}
e^{-\frac{\dim\lam}{4}\tr(x^2)}
e^{\dim\lam \sum_j\frac{t_j}{2j}\tr(x^j)}d\mu_\lam(x)\\
&\quad+\sum_{\lam\in\Ghat_4}
\log \int_{\mathcal{H}_{\dim\lam/2, \mathbb{H}}}
e^{-\frac{\dim\lam}{2}\tr(x^2)}
e^{\dim\lam \sum_j\frac{t_j}{j}
\tr_{\dim\lam/2}(x^j)}d\mu_\lam(x)\\
&=\sum_{\substack{\Gamma  \text{ connected}\\
\text{M\"obius graph}}}\frac{1}{|\Aut_M(\Gamma)|}
\sum_{\lam\in\Ghat_1} \;(\dim\lam)^{\rchi(S_\Gamma)}
\prod_j t_j ^{v_j(\Gamma)}\\
&\quad+\sum_{\substack{\Gamma  \text{ connected orientable}\\
\text{M\"obius graph}}}\frac{1}{|\Aut_M(\Gamma)|}
\sum_{\lam\in\Ghat_2} \;(\dim\lam)^{\rchi(S_\Gamma)}
\prod_j t_j ^{v_j(\Gamma)}\\
&\quad+\sum_{\substack{\Gamma  \text{ connected}\\
\text{M\"obius graph}}}\frac{1}{|\Aut_M(\Gamma)|}
\sum_{\lam\in\Ghat_4} \;(-\dim\lam)^{\rchi(S_\Gamma)}
\prod_j t_j ^{v_j(\Gamma)}\\
&=\sum_{\substack{\Gamma  \text{ connected orientable}\\
\text{M\"obius graph}}}\frac{1}{|\Aut_M(\Gamma)|}
\sum_{\lam\in\Ghat} \;(\dim\lam)^{\rchi(S_\Gamma)}
\prod_j t_j ^{v_j(\Gamma)}\\
&\quad+\sum_{\substack{\Gamma  \text{ connected non-}\\
\text{orientable M\"obius graph}}}\frac{1}{|\Aut_M(\Gamma)|}
\sum_{\lam\in\Ghat_1} \;(\dim\lam)^{\rchi(S_\Gamma)}
\prod_j t_j ^{v_j(\Gamma)}\\
&\quad+\sum_{\substack{\Gamma  \text{ connected non-}\\
\text{orientable M\"obius graph}}}\frac{1}{|\Aut_M(\Gamma)|}
\sum_{\lam\in\Ghat_4} \;(-\dim\lam)^{\rchi(S_\Gamma)}
\prod_j t_j ^{v_j(\Gamma)}\;.
\end{split}
\end{equation}
Notice that the sum over orientable M\"obius graphs
recovers Mednykh's formula (\ref{eq:Med}) again,
because the Euler characteristic $\rchi(S)$ is
even for an orientable surface. From the sum over
non-orientable M\"obius graphs, we obtain the formula
of Frobenius-Schur (\ref{eq:FScount}) of \cite{FS}:
\begin{equation*}
\sum_{\lam\in\Ghat_1}
(\dim\lam)^{\rchi(S)} 
+
\sum_{\lam\in\Ghat_4}
(-\dim\lam)^{\rchi(S)}
= |G|^{\rchi(S)-1}|\Hom(\pi_1(S),G)|   \;.
\end{equation*}
Note that the $\hat{G}_2$ component has no contribution in
this formula. This is due to the fact that graphical
expansion of a complex
hermitian matrix integral contains only orientable
 ribbon graphs.


\providecommand{\bysame}{\leavevmode\hbox to3em{\hrulefill}\thinspace}

\bibliographystyle{amsplain}

\begin{thebibliography}{10}




\bibitem{AV}
Mark Adler and Pierre van Moerbeke, 
\emph{Hermitian, symmetric and symplectic 
random ensembles: PDEs for the distribution of the spectrum},
Annals  of Mathematics (2) \textbf{153} (2001), no. 1, 149--189.

\bibitem{BDJ}
Jinho Baik, Percy Deift and Kurt Johansson,
\emph{On the
distribution of the length of the longest increasing
subsequence of random permutations}, Journal of  American
Mathematical  Society \textbf{ 12} (1999), 1119--1178.

\bibitem{BDJ2}
Jinho Baik, Percy Deift and Kurt Johansson,
\emph{
On the distribution of the length of the
second row of a Young diagram under Plancherel measure},
math.CO/9901118 (1999). 

\bibitem{BR}
Jinho Baik and Eric Rains, 
\emph{Symmetrized random permutations},
math.CO/9910019 (1999).

\bibitem{Belyi}
G.~V.~Belyi, \emph{On galois extensions of a maximal cyclotomic fields}, Math.\
  U.S.S.R.\ Izvestija \textbf{14} (1980), 247--256.

\bibitem{BIZ}
D.~Bessis, C.~Itzykson and J.~B.~Zuber,
\emph{Quantum field theory techniques in graphical
enumeration}, Advanced in Applied Mathematics
\textbf{1} (1980), 109--157.

\bibitem{BI}
Pavel M.~Bleher and Alexander R.~Its,
\emph{Random matrix models and their applications},
Mathematical Sciences Research Institute Publications
\textbf{40}, Cambridge University Press, 2001.



\bibitem{BOO}
Alexander Borodin, Andrei Okounkov, and  Grigori Olshanski, 
\emph{On
asymptotics of Plancherel measures for symmetric
groups}, math.CO/990532 (1999).



\bibitem{BIPZ}
C.~Br\'ezin, C.~Itzykson, G.~Parisi, 
and J.-B.~Zuber,
\emph{Planar diagrams},
Communications in Mathematical Physics \textbf{59} (1978), 35--51.






\bibitem{Burnside}
William Burnside, 
\emph{Theory of groups of finite order}, Second Edition, 
Cambridge University Press, 1991.




\bibitem{D}
Percy Deift, 
{Integrable systems and combinatorial theory},
Notices of AMS, \textbf{47} (2000), 631--640.

\bibitem{Feynman}
Richard P.~Feynman,
\emph{Space-time approach to quantum electrodynamics},
Physical Review \textbf{76} (1949), 769--789. 



\bibitem{FQ}
Daniel S. Freed and Frank Quinn, 
\emph{Chern-Simons theory with finite gauge group},
Communications in Mathematical Physics \textbf{156}
(1993), 435--472.



\bibitem{F}
Georg Frobenius, 
\emph{\"Uber  {G}ruppencharaktere},
Sitzungsberichte der k\"oniglich preussischen 
{A}kademie der {W}issenschaften 
(1896),  985--1021.

\bibitem{FS}
Georg Frobenius and Isaai Schur, 
\emph{\"Uber die reellen {D}arstellungen der 
endlichen {G}ruppen},
Sitzungsberichte der k\"oniglich preussischen 
{A}kademie der {W}issenschaften 
(1906), 186--208.

\bibitem{Gardiner}
Frederick~P. Gardiner, \emph{{Teichm\"uller} 
theory and quadratic
  differentials}, John Wiley \& Sons, 1987.



\bibitem{GHJ}
I.~P.~Goulden,  J.~L.~Harer, J.~L. and D.~M.~Jackson, 
\emph{A geometric parametrization for the virtual {E}uler
characteristics of the moduli spaces of real and complex
algebraic curves},
Trans.\ Amer.\ Math.\ Soc.\
\textbf{353} (2001), 4405--4427.

\bibitem{GrossTaylor}
David~J.~Gross and Washington Taylor, 
\emph{Two dimensional QCD is a string theory},
Nucl.\ Phys.\ \textbf{B 400} (1993), 181--210.


\bibitem{GT}
Jonathan L.~Gross and Thomas W.~Tucker,
\emph{Topological graph theory}, John Wiley \& Sons, 1987.


\bibitem{Grothendieck}
Alexander Grothendieck, 
\emph{Esquisse d'un programme} (1984), 
reprinted in \cite {SL}, 7--48.

\bibitem{Harer1986}
John~L.~Harer, \emph{The virtual cohomological
dimension of the mapping class group of an
orientable surface},
     Inventiones Mathematicae \textbf{84} (1986),
157--176.






\bibitem{Harer-Zagier}
John~L.\ Harer and Don Zagier, 
\emph{The {Euler} 
characteristic of the moduli
  space of curves}, 
Inventiones Mathematicae 
\textbf{85} (1986), 457--485.

\bibitem{Hatcher}
Allen  Hatcher,
\emph{On triangulations of surfaces},
Topology and it Applications \textbf{40} (1991),
189--194. 


\bibitem{Isaacs}
I.~Martin Isaacs, \emph{Character theory of finite groups},
Academic Press, 1976.

\bibitem{J}
Kurt Johansson, 
\emph{Discrete orthogonal polynomial
ensembles and the Plancherel measure}, math.CO/9906120 (1999). 

\bibitem{Jones}
Gareth A.~Jones, \emph{Characters and surfaces: a survey},
London Mathematical Society Lecture Note 
Series \textbf{249}, 
\emph{The atlas of finite groups: ten years on},
Robert Curtis and Robert Wilson, Eds., (1998), 90--118.



\bibitem{Kontsevich}
Maxim Kontsevich, 
\emph{Intersection theory on 
the moduli space of 
curves and the
matrix {Airy} function},
  Communications in Mathematical Physics
 \textbf{147} (1992), 1--23.

\bibitem{L1}
Kefeng Liu, 
\emph{Heat kernel and moduli space},
Mathematical Research Letters \textbf{3}
(1996), 743--762.

\bibitem{L2}
Kefeng Liu, 
\emph{Heat kernel and moduli space II},
Mathematical Research Letters \textbf{4}
(1996), 569--588.



\bibitem{L3}
Kefeng Liu, 
\emph{Heat kernels, symplectic geometry, 
moduli spaces and finite groups},
Surveys in Differential Geometry \textbf{5}
(1999), 527--542.




\bibitem{Med}
A.~D.~Mednykh, 
\emph{Determination of the number of 
nonequivalent coverings over a compact Riemann surface},
Soviet Mathematics Doklady \textbf{19} (1978), 318--320.


\bibitem{Mehta}
Madan Lal Mehta, 
\emph{Random matrices}, Second Edition, 
Academic Press, 1991.


\bibitem{Mulase1994a} 
Motohico Mulase, \emph{
Algebraic theory of the KP equations},
in Perspectives in Mathematical Physics, 
R.~Penner and S.~T.~Yau, Editors., 
Intern.\ Press Co. (1994),  157--223.  

\bibitem{Mulase1994b} 
Motohico Mulase, \emph{Matrix integrals
 and integrable systems},
  in Topology, geometry and field theory, K.~Fukaya
et al.\ Editors,
World Scientific (1994), 111--127.




\bibitem{Mulase1995}
Motohico Mulase, 
\emph{Asymptotic analysis of a 
hermitian matrix integral},
  International Journal of Mathematics 
\textbf{6} (1995), 881--892.

\bibitem{Mulase1998}
Motohico Mulase, \emph{
Lectures on the asymptotic
 expansion of a hermitian matrix integral},
in Supersymmetry and Integrable Models, 
Henrik Aratin et al., Editors, 
Springer Lecture Notes in Physics
\textbf{502}  (1998), 91--134.


\bibitem{MP1998}
Motohico Mulase and Michael Penkava,
\emph{Ribbon graphs, quadratic differentials on
{R}iemann surfaces, and algebraic curves defined
over $\overline{\mathbb{Q}}$}, Asian Journal of
Mathematics \textbf{2} (1998), 875--920.



\bibitem{MW}
Motohico Mulase and Andrew Waldron,
\emph{Duality of orthogonal and symplectic matrix 
integrals and quaternionic {F}eynman graphs}, 
math-ph/0206011 (2002).

\bibitem{MY}
Motohico Mulase and Josephine T.~Yu,
\emph{A 
generating function of 
the number of homomorphisms from a 
surface group into a finite group}, 
math.QA/0209008 (2002).



\bibitem{O}
Andrei Okounkov,
\emph{Random matrices and random permutations},
math.CO/9903176 (1999). 



\bibitem{OP1}
Andrei Okounkov and Rahul Pandharipande,
\emph{Gromov-Witten theory, Hurwitz numbers, and matrix 
models, I}, math.AG/0101147 (2001).



\bibitem{OP2}
Andrei Okounkov and Rahul Pandharipande,
\emph{The equivariant Gromov-Witten theory of 
$\mathbb{P}^1$}, 
math.AG/0207233 (2002).



\bibitem{Penner}
Robert~C.\  Penner, 
\emph{Perturbation series and the moduli 
space of {Riemann}
  surfaces}, Journal of Differential 
Geometry \textbf{27}
 (1988), 35--53.



\bibitem{Ringel}
Gerhard Ringel, \emph{Map color theorem},
Springer-Verlag, 1974.



\bibitem{Schneps}
Leila Schneps, 
\emph{The grothendieck theory of dessins d'enfants}, 
  London Mathematical Society Lecture Notes Series,
vol.\ 200, 1994.

\bibitem{SL}
Leila Schneps and Pierre Lochak, editors,
\emph{Geometric Galois actions: Around
Grothendieck's esquisse d'un programme},
      London Mathematical Society Lecture Notes Series,
vol.\ 242, 1997.




\bibitem{Serre}
Jean-Pierre Serre,
\emph{Linear representations of finite groups},
Springer-Verlag, 1987.

\bibitem{Stanley}
Richard P.\ Stanley, 
\emph{Enumerative combinatorics}
volume 2, Cambridge University Press, 2001.

\bibitem{Strebel}
Kurt Strebel, 
\emph{Quadratic differentials}, Springer-Verlag, 1984.


\bibitem{'tHooft}
Gerard 't~Hooft,
\emph{A planer diagram theory for
strong interactions},
Nuclear Physics B
\textbf{72} (1974),
461--473.

\bibitem{TW}
Craig A.~Tracy and Harold Widom, 
\emph{Fredholm Determinants, Differential Equations and Matrix Models},
hep-th/9306042 Communications in Mathematical
Physics \textbf{163} (1994), 33--72.


\bibitem{VM}
Pierre van Moerbeke, 
\emph{Integrable lattices: random matrics and random permutations},
in Random Matrix Models and Their Applications, Bleher and Its,
Editors, MSRI Publications \textbf{40} (2001), 321--406.



\bibitem{W1991}
Edward Witten, \emph{On quantum gauge theories in
two dimensions}, 
 Communications in Mathematical Physics \textbf{141} (1991),
153--209.


\bibitem{W1991a}
Edward Witten, \emph{Two dimensional gravity and
intersection
theory on moduli space}, Surveys in
Differential Geometry \textbf{1} (1991),
243--310.

\bibitem{Yu}
Josephine Yu, 
\emph{Graphical expansion of matrix integrals with
values in a Clifford algebra}, in Explorations: 
A Journal of Undergraduate Research, 
University of California, Davis,  vol.~\textbf{6}
(2003).



\end{thebibliography}

\end{document}